\documentclass[10pt,twoside,a4paper,reqno]{amsart}
\usepackage{amscd,amsmath,amsthm,amsfonts,latexsym,amssymb}

\theoremstyle{plain}
\newtheorem{theo+}           {Theorem}
\newtheorem{prop+}           {Proposition}
\newtheorem{coro+}           {Corollary}
\newtheorem{lemm+}           {Lemma}
\newtheorem{conjecture}      {Conjecture}

\theoremstyle{definition}
\newtheorem{defi+}           {Definition}

\newtheorem{not+}            {Notation}
\newtheorem{ex}              {Example}

\theoremstyle{remark}
\newtheorem{rema+}           {Remark}

\newenvironment{theorem}{\begin{theo+}}{\end{theo+}}

\newenvironment{corollary}{\begin{coro+}}{\end{coro+}}
\newenvironment{lemma}{\begin{lemm+}}{\end{lemm+}}
\newenvironment{remark}{\begin{rema+}}{\end{rema+}}
\newenvironment{definition}{\begin{defi+}}{\end{defi+}}
\newenvironment{notation}{\begin{not+}}{\end{not+}}

\newcommand {\bC} {\mathbb {C}}
\newcommand {\bD} {\mathbb {D}}
\newcommand {\bR} {\mathbb {R}}

\newcommand {\bX} {\mathbb {X}}
\newcommand {\bN} {\mathbb {N}}
\newcommand {\bM} {\mathbb {M}}

\newcommand {\al} {\alpha}
\newcommand {\be} {\beta}
\newcommand {\ga} {\gamma}
\newcommand {\de} {\delta}
\newcommand {\la} {\lambda}
\newcommand {\ze} {\zeta}
\newcommand {\si} {\sigma}
\newcommand {\om} {\omega}
\newcommand {\te} {\theta}

\newcommand {\Om} {\Omega}
\newcommand {\De} {\Delta}
\newcommand {\Si} {\Sigma}

\newcommand {\calH} {\mathcal {H}}

\newcommand {\calM} {\mathcal {M}}
\newcommand {\calN} {\mathcal {N}}
\newcommand {\calK} {\mathcal {K}}
\newcommand {\calL} {\mathcal {L}}
\newcommand {\calZ} {\mathcal {Z}}
\newcommand {\calA} {\mathcal {A}}
\newcommand {\calX} {\mathcal {X}}
\newcommand {\calO} {\mathcal {O}}
\newcommand {\calW} {\mathcal {W}}

\newcommand{\co} {\text{conv}}
\newcommand{\tco} {\text{{\em conv}}}

\newcommand{\ba}{\mathbf a}
\newcommand{\bb}{\mathbf b}
\newcommand{\bx}{\mathbf x}
\newcommand{\by}{\mathbf y}
\newcommand{\bu}{\mathbf u}
\newcommand{\bv}{\mathbf v}

\newcommand{\bs}{\mathbf s}
\newcommand{\bi}{\mathbf i}
\newcommand{\bj}{\mathbf j}
\newcommand{\bbe}{\mathbf e}
\newcommand{\bbf}{\mathbf f}
\newcommand{\bz}{\mathbf z}
\newcommand{\pr}{\prec}

\begin{document}

\numberwithin{equation}{section}

\title[Equilibrium points of logarithmic potentials. I.]
{Equilibrium points of logarithmic potentials induced by positive charge
distributions. I. Generalized de Bruijn-Springer relations}
\author[J.~Borcea]{Julius Borcea}
\address{Department of Mathematics, Stockholm University, SE-106 91 Stockholm,
   Sweden}
\email{julius@math.su.se}


\begin{abstract}
A notion of weighted multivariate majorization is defined as a preorder on 
sequences of vectors in Euclidean space induced by the Choquet ordering for 
atomic probability measures. We characterize this preorder both in terms of 
stochastic matrices and convex functions and use it to describe the 
distribution of equilibrium points of logarithmic potentials generated by 
discrete planar charge configurations. In the case of $n$ positive charges we 
prove that the equilibrium points satisfy $\binom{n}{2}$ weighted 
majorization relations and are uniquely determined by $n-1$ such relations. 
It is further shown that the Hausdorff geometry of the equilibrium points and 
the charged particles is controlled by the weighted standard deviation of the 
latter. By using finite-rank perturbations of compact normal Hilbert space 
operators we establish similar relations for infinite charge distributions. 
We also discuss a hierarchy of weighted de Bruijn-Springer relations and 
inertia laws, the existence of zeros of Borel series with positive 
$l^1$-coefficients, and an operator version of the Clunie-Eremenko-Rossi 
conjecture. 
\end{abstract}

\maketitle

\tableofcontents

\section{Introduction}\label{s0}

Let $\{z_i\}_{i\in\calN}$ be a sequence of distinct points in a domain 
$\Omega\subseteq \bC$, where $\calN\subseteq \bN$ with $|\calN|\ge 2$. If at 
each $z_i$ we 
place a particle with charge $a_i$ -- actually, an infinitely long thin rod 
perpendicular to the plane with uniform charge distribution -- then up to
complex conjugation and a constant factor the resulting electrostatic force is 
given by 
\begin{equation}\label{n-f}
f(z)=\sum_{i\in\calN}\frac{a_i}{z-z_i}
\end{equation}
and the logarithmic potential generated by these charged particles is 
\begin{equation}\label{n-pot}
U(z)=\sum_{i\in\calN}a_{i}\log\left|1-\frac{z}{z_i}\right|
\end{equation}
(cf.~\cite[p.~10]{Ke}). If $|\calN|=\infty$ then the series in~\eqref{n-f} 
and~\eqref{n-pot} converge absolutely for all 
$z\in\Omega\setminus\{z_i\}_{i\in\calN}$ provided that certain 
growth conditions are satisfied (see \S \ref{s22} below). In this 
case $f$ is meromorphic while $U$ is subharmonic in $\Omega$. Since 
$\nabla U(z)=\overline{f(z)}$ the critical points of the 
potential coincide with the equilibrium points of the field, 
that is, points upon which a free electron (or rod) once placed would remain.
An important problem that goes back to Gauss is to describe 
the distribution 
and the relative geometry of the equilibrium points and charged particles. One 
of the main purposes of this paper is to study this problem by means of the 
{\em Hausdorff distance} and the notion of {\em weighted multivariate 
majorization} that we define in \S \ref{s1} as an extension
of classical and matrix majorization. It is interesting to note that 
although our results are purely analytic or geometric in nature, they 
are obtained through a combination of the multilinear algebra and 
operator theoretic methods developed in \cite{Ma1,Ma2} and \cite{Ma1,Ma2,P},
respectively. The fact that such geometric properties could be established 
in this way was actually first noticed in {\em loc.~cit.}, where 
long-standing conjectures on the geometry of polynomials were solved by 
applying this type of techniques. Below we make extensive use of the key 
observation made in the aforementioned papers, namely that studying 
zeros and critical points of complex polynomials -- or, more generally, 
poles and zeros of logarithmic potentials generated by finite positive charge 
configurations -- is equivalent to studying spectra of normal matrices and 
their compressions.

The weighted multivariate generalization of vector and matrix majorization 
introduced in \S \ref{s1} is a natural preorder on finite sequences of vectors
in Euclidean space induced by the Choquet ordering for nonnegative Radon 
measures (Definition~\ref{def-w}). We characterize this 
new notion in terms of stochastic matrices and convex functions and we 
recover in this way several related results on classical and multivariate 
majorization. 
In the process we use a Hahn-Banach type argument to give a new simple proof 
of a fundamental theorem of Sherman (Theorem~\ref{gen-S}).

In \S \ref{s2} we state most of our main results. These describe the 
distribution and the relative location of equilibrium points of logarithmic 
potentials for discrete bounded planar configurations of positive 
charges and are of three different types. First, we consider $n$ positive 
charges and prove that the equilibrium points satisfy $\binom{n}{2}$ weighted 
multivariate majorization relations and are in fact uniquely determined 
by $n-1$ such relations (Theorems~\ref{t-1}--\ref{t-4}). The higher 
de Bruijn-Springer relations thus obtained generalize some of the results of
\cite{Ma1,Ma2} to arbitrary (not necessarily equal) positive charges. 
Second, we show that 
the relative Hausdorff geometry of the extended set of equilibrium points and 
the set of charged particles is controlled by the weighted standard deviation 
of the latter (Theorems~\ref{t-5}--\ref{t-6}). For equal charges the bounds 
given by these theorems 
are much better and also more natural than those predicted by the famous  
Sendov conjecture for complex polynomials and its converse, 
see, e.g., \cite{Bo2}. As we explain in the sequel \cite{BP2} to this 
paper, these results also seem to suggest that the Hausdorff geometry is 
governed by natural dispersion measures associated with the given charge 
configuration. Finally, we establish weighted 
multivariate majorization relations for essentially any bounded discrete set 
of positively charged particles with finite total charge for which 
equilibrium points are known to exist (Theorems~\ref{t-7}--\ref{t-8}). We 
further conjecture that equilibrium points actually always exist under some 
mild assumptions (Conjectures~\ref{con-2}--\ref{new-c2} in \S \ref{s42}). 

The proofs of the aforementioned results are given in \S \ref{s3}. As we 
already pointed out, for finite charge configurations these are based on the 
multilinear algebra and matrix theory techniques of \cite{Ma1,Ma2,P}. In the 
process we are led to consider several natural related questions, such as a 
hierarchy of inertia laws for finite systems of planar charged particles
that we propose and discuss in \S \ref{s-mom} (Conjecture~\ref{gen-mom}). 
This would yield a whole new series of weighted generalized 
de Bruijn-Springer relations. The techniques mentioned above actually extend 
to the case of infinite discrete distributions of charge, which allows us to
make use of operator theoretic methods such as 
compressions of compact normal operators on separable Hilbert space in order 
to deal with this case. In this 
way we arrive at a great many related questions on the interplay between value 
distribution of meromorphic functions and variations of discrete spectra of 
compact normal Hilbert space operators under finite-rank perturbations. In 
\S \ref{s42}--\ref{s43} we give an operator theoretic approach to  
the well-known Clunie-Eremenko-Rossi conjecture \cite{CER} dealing with the 
existence of zeros of certain Borel series and we consider possible extensions 
of this conjecture as well as their potential implications 
(Conjectures~\ref{con-CER}--\ref{con-compl}).

\section{Weighted multivariate majorization}\label{s1}

The concept of majorization was first studied by economists 
early in the twentieth century as a means for altering the unevenness of 
distribution of wealth or income. The majorization preorder on $n$-tuples of
real numbers -- also known as the strong spectral order, vector majorization 
or classical majorization -- was used 
by Hardy, Littlewood and P\'olya in their study of analytic inequalities 
\cite{HLP} and by Schur in his work on the spectra of positive semidefinite 
Hermitian operators (cf.~\cite{DK}). Classical majorization has since become 
an important tool in analysis, combinatorics and 
statistics \cite{MO} and has recently found remarkable applications 
to quantum mechanics and entanglement transformations in quantum 
computation and information theory (see \cite{Bo1} and references therein).
The strong spectral order essentially quantifies the intuitive notion that
the components of a real $n$-vector $\bx$ are ``less spread out'' than the
components of another such vector $\by$. Several matrix versions of this 
notion have been proposed and studied in various contexts 
(cf., e.g., \cite{Da} and \cite{MO}). 

In this section we generalize classical and matrix 
majorization by introducing a notion of 
weighted multivariate majorization that we shall use throughout both this
paper and its sequel \cite{BP2}.

\begin{notation}\label{new-not1}
Let $k\in \bN$. For each $n\in \bN$ we define
{\allowdisplaybreaks
\begin{align}\label{sp-n}
&\calX_{n}^{k}=\left\{(\bx_1,\ldots,\bx_n)\,\,\big|\,\,\bx_i\in \bR^k,
1\le i\le n\right\}\\
&\calA_{n}=\left\{(a_1,\ldots,a_n)\,\,\bigg|\,\,a_i\in (0,1),1\le i\le n,
\sum_{i=1}^{n}a_i=1\right\}\notag\\
&\bX_{n}^{k}=\calX_{n}^{k}\times \calA_{n}\notag
\end{align}
and set
\begin{equation}\label{sp}
\bX^{k}=\bigcup_{n=1}^{\infty}\bX_{n}^{k}.
\end{equation}
Given $X=(\bx_1,\ldots,\bx_m)\in \calX_{m}^k$ we denote by 
$X^T=(\bx_1,\ldots,\bx_m)^T$ its ``transpose'', i.e., the ordered 
$m$-tuple consisting of the vectors $\bx_1,\ldots,\bx_m$ written in column 
form.} We identify $X^T$ with the $m\times k$ real matrix 
whose $i$-th row consists of the coordinates of the vector $\bx_i$ in the 
standard basis of $\bR^k$. Let finally 
$\bM_{m,n}^{rs}$ be the set of all row stochastic $m\times n$ matrices.
\end{notation} 

The notion of weighted multivariate majorization defined below 
(see Definition~\ref{def-w}) is motivated 
by the following theorem of Sherman \cite{Sh}:

\begin{theorem}\label{gen-S}
Let $(X,\ba)\in \bX_{m}^k$ and $(Y,\bb)\in \bX_{n}^k$, where 
$X=(\bx_1,\ldots,\bx_m)\in \calX_{m}^k$, $\ba=(a_1,\ldots,a_m)\in \calA_m$,
$Y=(\by_1,\ldots,\by_n)\in \calX_{n}^k$ and $\bb=(b_1,\ldots,b_n)\in \calA_n$.
The following conditions are equivalent:
\begin{itemize}
\item[(i)] For any (continuous) convex function $\Phi:\bR^k\rightarrow \bR$ 
one has 
$$\sum_{i=1}^{m}a_{i}\Phi(\bx_i)\le \sum_{j=1}^{n}b_{j}\Phi(\by_j).$$
\item[(ii)] There exists a matrix $R\in \bM_{m,n}^{rs}$ 
such that 
$$\tilde{X}^T=R\tilde{Y}^T\text{ and }\,\bb=\ba R,$$
where $\tilde{X}^T$ is an $m\times k$ matrix 
and $\tilde{Y}^T$ is an $n\times k$ matrix obtained by some (and then any) 
ordering of the vectors in $X^T$ and $Y^T$, respectively.
\end{itemize}
\end{theorem}

The special case when $m=n$ and all weights are 
equal in Theorem~\ref{gen-S} was long assumed to be an open question and 
appears only as an implicit conjecture in \cite[p. 433]{MO}, which is the 
definite reference on majorization theory. Although essentially correct, 
Sherman's paper used rather complicated arguments and contained a series of 
misprints that required several subsequent errata, which may explain why the 
status of his proof seems to have remained unclear for such a long period. 
Below we present a short and simple proof of Theorem~\ref{gen-S}.

Before embarking on the actual proof let us note that condition (i) in 
Sherman's theorem may also be written in terms of the Choquet ordering for 
positive measures, a notion that we shall actually use in~\cite{BP2} (see
Definition 9 in {\em loc.~cit.}). Indeed, any pair $(X,\ba)\in\bX_{m}^{k}$ 
uniquely determines a discrete probability measure $\mu_{_{(X,\ba)}}$ on 
$\bR^k$ concentrating a mass $a_i$ at $\bx_i$, $1\le i\le m$. Given two 
positive measures $\mu_1$ and $\mu_2$ on $\bR^k$ one says that $\mu_2$ 
{\em dominates} $\mu_1$ {\em in the Choquet ordering} or that $\mu_2$ is a 
{\em dilation} of $\mu_1$, denoted $\mu_1\pr \mu_2$, 
if $\mu_1(\Phi)\le \mu_2(\Phi)$ for any continuous convex function $\Phi$ on 
$\bR^k$. Theorem~\ref{gen-S} (i) may therefore be rewritten as 
$\mu_{_{(X,\ba)}}\pr \mu_{_{(Y,\bb)}}$. We refer to~\cite{BR,CFM,EH,FH,Ph} 
for a detailed account of the Choquet ordering and various comparisons of 
nonnegative Radon measures and probability distributions in locally convex 
separable topological spaces. In particular, in the latter setting 
Sherman's theorem may in fact be viewed as a special case of the 
Cartier-Fell-Meyer theorem \cite{CFM,Ph,BR} and its subsequent extension due 
to Fischer-Holbrook \cite{FH}. The most general result of this type, namely a 
basic representation theorem for convex domination of measures, seems to be 
contained in \cite{EH}.

\begin{proof}[Proof of Theorem~\ref{gen-S}]
Since the implication (ii) $\Rightarrow$ (i) is an immediate consequence of 
Jensen's inequality we shall only focus on the converse statement. Assume that
condition (i) in Theorem~\ref{gen-S} holds, 
let $X=(\bx_1,\ldots,\bx_m)\in \calX_{m}^k$, $\ba=(a_1,\ldots,a_m)\in\calA_m$, 
$Y=(\by_1,\ldots,\by_n)\in \calX_{n}^k$, $\bb=(b_1,\ldots,b_n)\in \calA_n$,  
and define the following set
$$\calM(X,Y)=\left\{A=(a_{ij})\in \bM_{m,n}^{rs}\,\,\bigg|\,\,
\bx_i=\sum_{j=1}^{n}a_{ij}\by_j,\,1\le i\le m\right\}.$$

\begin{lemma}\label{l1}
$\calM(X,Y)$ is a non-empty closed convex subset of $\bM_{m,n}^{rs}$.
\end{lemma}

\begin{proof}
It is clear that $\calM(X,Y)$ is closed and convex. To show that 
$\calM(X,Y)\neq \emptyset$ it is enough to prove that if 
Theorem~\ref{gen-S} (i) holds then $\bx_i\in \co(\by_1,\ldots,\by_n)$ for 
$1\le i\le m$, where $\co(\by_1,\ldots,\by_n)$
is the convex hull of the vectors $\by_1,\ldots,\by_n$. Supposing that 
$\bx_l\notin \co(\by_1,\ldots,\by_n)$ for some $l\in \{1,\ldots,m\}$ let $d$ 
denote the Euclidean distance in $\bR^k$ and consider the function
$\Phi:\bR^k\rightarrow \bR$ given by 
$\Phi(\bu)=d(\bu,\co(\by_1,\ldots,\by_n))$. 
Note that $\Phi(\bx_l)>0$, $\Phi(\by_j)=0$ for $1\le j\le n$, and $\Phi$ is
a convex function since $\co(\by_1,\ldots,\by_n)$ is a non-empty convex subset 
of $\bR^k$. It follows that 
$\sum_{i=1}^{m}a_{i}\Phi(\bx_i)>\sum_{j=1}^{n}b_{j}\Phi(\by_j)$,
which contradicts condition (i) of Theorem~\ref{gen-S}. Thus
 $\bx_i\in \co(\by_1,\ldots,\by_n)$ for $1\le i\le m$, as required.
\end{proof}

Given $M=(m_{ij})\in \calM(X,Y)$ we define a vector 
$\bs(M)=(s_1(M),\ldots,s_n(M))\in \bR^n$ by 
$s_j(M)=\sum_{i=1}^{m}a_im_{ij}$, $1\le j\le n$. Set
$$\calL(X,Y)=\{\bs(M)\mid M\in\calM(X,Y)\}.$$

Note that $\al\bs(M)+(1-\al)\bs(N)=\bs(\al M+(1-\al)N)$ for any 
$M,N\in \calM(X,Y)$ and $\al\in [0,1]$, so that $\calL(X,Y)$ is a convex 
subset of $\bR^n$. Moreover, since $\calM(X,Y)$ is closed the same must be 
true for $\calL(X,Y)$. 

\begin{lemma}\label{l2}
With the above notations one has $\bb\in \calL(X,Y)$. In other words, 
there exists a matrix $M_0\in \calM(X,Y)$ such that $\bs(M_0)=\bb$.
\end{lemma}

\begin{proof}
Assume that $\bb\notin \calL(X,Y)$. Since $\calL(X,Y)$ is a closed convex set 
it follows from the Hahn-Banach theorem that there exist $c\in \bR$ and 
$(r_1,\ldots, r_n)\in \bR^n$ such that
\begin{equation}\label{h-b}
\sum_{j=1}^{n}r_js_j<c<\sum_{j=1}^{n}r_{j}b_j\text{ for }
(s_1,\ldots,s_n)\in \calL(X,Y).
\end{equation}
Given an arbitrary vector 
$\by\in \co(\by_1,\ldots,\by_n)$ we define the simplex
$$W(\by)=\left\{(\la_1,\ldots,\la_n)\in [0,1]^n\,\,\bigg|\,\,
\by=\sum_{j=1}^{n}\la_j
\by_j,\sum_{j=1}^{n}\la_j=1\right\}.$$
Note that $W(\by)$ is compact for any $\by\in \co(\by_1,\ldots,\by_n)$
since it is complete and totally bounded. We may therefore define a function 
$G:\co(\by_1,\ldots,\by_n)\rightarrow \bR$ by setting
$$G(\by)=\max\left\{\sum_{j=1}^{n}r_j\la_j\,\,\bigg|\,\,(\la_1,\ldots,\la_n)
\in W(\by)\right\}\text{ for }\by\in \co(\by_1,\ldots,\by_n).$$
Let $\la_j(\by)\in [0,1]$, $1\le j\le n$, be such that 
$G(\by)=\sum_{j=1}^{n}r_j\la_j(\by)$. Since $\la_j(\by_i)=\de_{ij}$ one has
$G(\by_j)\ge r_j$ for $1\le j\le n$ and thus
\begin{equation}\label{h-b1}
\sum_{j=1}^{n}b_{j}G(\by_j)\ge \sum_{j=1}^{n}b_{j}r_{j}>c.
\end{equation}
On the other hand the left inequality in~\eqref{h-b} implies that for any 
matrix $M=(m_{ij})\in \calM(X,Y)$ one has
\begin{equation}\label{h-b2}
c>\sum_{j=1}^{n}r_js_j(M)=\sum_{i=1}^{m}\sum_{j=1}^{n}a_{i}r_{j}m_{ij}.
\end{equation}
Let us now consider the extension $\tilde{G}$ of the function $G$ to $\bR^k$
given by $\tilde{G}(\by)=G(\by)$ if $\by\in \co(\by_1,\ldots,\by_n)$ and 
$\tilde{G}(\by)=0$
for $\by\in \bR^k\setminus \co(\by_1,\ldots,\by_n)$. It is not difficult to 
see that 
$\tilde{G}$ is a concave function on $\bR^k$. By Lemma~\ref{l1} one has 
$\bx_i\in \co(\by_1,\ldots,\by_n)$, $1\le i\le m$, so that the numbers 
$\la_j(\bx_i)$ are well defined for $1\le i\le m$, $1\le j\le n$. Note that
$(\la_j(\bx_i))\in \calM(X,Y)$ and that by~\eqref{h-b2} and~\eqref{h-b1} the 
following holds
$$\sum_{i=1}^{m}a_{i}\tilde{G}(\bx_i)
=\sum_{i=1}^{m}\sum_{j=1}^{n}a_{i}r_j\la_j(\bx_i)
<\sum_{j=1}^{n}b_{j}r_{j}\le \sum_{j=1}^{n}b_{j}\tilde{G}(\by_j).$$
This contradicts however condition (i) of Theorem~\ref{gen-S} and thus we are 
done. It follows that
$\bb\in \calL(X,Y)$, which completes the proof of the lemma.
\end{proof}
The fact that (ii) $\Rightarrow$ (i) in Theorem~\ref{gen-S} is equivalent 
to Lemma~\ref{l2}.
\end{proof}

\begin{remark}\label{r-thanks}
We are grateful to the referee for pointing out the similarity 
between 
the Hahn-Banach type argument used above and the one given in the original 
proof of the Cartier-Fell-Meyer theorem \cite{CFM}.
\end{remark}

\begin{definition}\label{def-w}
The pair $(X,\ba)\in \bX_{m}^k$ is said to be {\em weightily majorized} by 
the pair $(Y,\bb)\in \bX_{n}^k$, denoted $(X,\ba)\pr (Y,\bb)$, if the 
conditions of Theorem~\ref{gen-S} are satisfied.
\end{definition}

\begin{remark}\label{bary}
If $(X,\ba)\pr (Y,\bb)$ then the $\ba$-barycenter of $X$ must coincide with the
$\bb$-barycenter of $Y$, that is, 
$\sum_{i=1}^{m}a_{i}\bx_i=\sum_{j=1}^{n}b_{j}\by_j$.
\end{remark}

It is clear from Definition~\ref{def-w} that the weighted majorization 
relation is a preorder on $\bX^k$. Indeed, by Theorem~\ref{gen-S} 
this relation is both reflexive and transitive. Moreover, if $m=n$ and 
$a_i=b_j=\frac{1}{n}$, $1\le i,j\le n$, then 
Theorem~\ref{gen-S} and Birkhoff's theorem \cite[Theorem A.2]{MO} imply that 
the weighted multivariate majorization relation induces a partial ordering on 
the set of unordered $n$-tuples of vectors in $\bR^k$. We shall refer to this 
partial ordering as the {\em ordinary multivariate majorization} relation, for 
which Sherman's theorem takes a particularly simple form:

\begin{corollary}\label{cor-po}
If $X=(\bx_1,\ldots,\bx_n)\in \calX_{n}^{k}$ and $Y=(\by_1,\ldots,\by_n)\in 
\calX_{n}^{k}$ then the following conditions are equivalent:
\begin{itemize}
\item[(i)] The inequality $\sum_{i=1}^{n}\Phi(\bx_i)\le 
\sum_{j=1}^{n}\Phi(\by_j)$
holds for any convex function $\Phi:\bR^k\rightarrow \bR$.
\item[(ii)] There exists a doubly stochastic $n\times n$ matrix $S$ such that 
$\tilde{X}^T=S\tilde{Y}^T$, where $\tilde{X}^T$ and $\tilde{Y}^T$ are 
$n\times k$ matrices obtained by some ordering of the vectors in 
$X^T$ and $Y^T$, respectively.
\end{itemize}
\end{corollary}

\begin{remark}
For $k=1$ Corollary~\ref{cor-po} amounts to a well-known
description of classical (vector) majorization due to Schur and to 
Hardy-Littlewood-P\'olya \cite{HLP,MO}. 
\end{remark}

\begin{remark}\label{ordi}
In the sequel we shall often consider weighted pairs of the form 
$(X,\ba)$ and $(Y,\bb)$, where $X$ and $Y$ are unordered tuples of 
complex numbers. In this
case we define the weighted majorization relation $(X,\ba)\pr (Y,\bb)$, 
when appropriate, by identifying $\bC$ with $\bR^2$ in Theorem~\ref{gen-S}
and Definition~\ref{def-w}. 
\end{remark}

Another interesting consequence of Sherman's theorem is as follows:

\begin{corollary}\label{cor-gs}
Assume that $m<n\le k$ and let $(\bx_1,\ldots,\bx_m)\in 
\calX_{m}^k$, $(\by_1,\ldots,\by_n)\in \calX_{n}^k$, 
$(a_1,\ldots,a_n)\in \calA_n$ and $(b_1,\ldots,b_n)\in \calA_n$. For
$m+1\le i\le n$ let $(c_{i1},\ldots,c_{in})\in \calA_n$ be such that
$b_{j}':=b_j-\sum_{i=m+1}^{n}a_{i}c_{ij}\neq 0$, $1\le j\le n$, and set 
$\bx_i=\sum_{j=1}^{n}c_{ij}\by_j$. If $\by_1,\ldots,\by_n$ are linearly 
independent then the following conditions are equivalent:
\begin{itemize}
\item[(i)] $\sum_{i=1}^{n}a_{i}\Phi(\bx_i)\le \sum_{j=1}^{n}b_{j}\Phi(\by_j)$ 
for all convex functions $\Phi:\bR^k\rightarrow \bR$.
\item[(ii)] $\sum_{i=1}^{m}a_i\Phi(\bx_i)\le \sum_{j=1}^{n}b_{j}'\Phi(\by_j)$ 
for all convex functions $\Phi:\bR^k\rightarrow \bR$.
\end{itemize}
If either of these conditions is fulfilled then $b_{j}'>0$ for $1\le j\le n$.
\end{corollary}

\begin{proof}
The implication (ii) $\Rightarrow$ (i) being trivial let us show that
(i) $\Rightarrow$ (ii). Set $X=(\bx_1,\ldots,\bx_n)\in \calX_{n}^{k}$, 
$X'=(\bx_1,\ldots,\bx_m)$, $Y=(\by_1,\ldots,\by_n)$, $\ba=(a_1,\ldots,a_n)$
and $\bb=(b_1,\ldots,b_n)$. Then condition (i) is equivalent to 
$(X,\ba)\pr (Y,\bb)$. Note that 
$\sum_{j=1}^{n}b_{j}'=\sum_{i=1}^{m}a_i=:\al\in (0,1)$ 
and let $\ba'=(\al^{-1}a_1,\ldots,\al^{-1}a_m)\in\calA_m$ and 
$\bb'=(\al^{-1}b_{1}',\ldots,\al^{-1}b_{n}')$. Using Theorem~\ref{gen-S} and 
the fact that $\by_1,\ldots,\by_n$ are linearly independent one can easily
show that $\bb'\in \calA_n$ and $(X',\ba')\pr (Y,\bb')$, which
is equivalent to condition (ii). 
\end{proof}

\begin{remark}
For an arbitrarily given convex function $\Phi:\bR^k\rightarrow \bR$
the inequality in Corollary~\ref{cor-gs} (i) is {\em \`a priori} 
weaker than the one in (ii). 
As a special case let us consider three complex 
numbers $z_1,z_2,z_3$ with $z_2\bar{z}_3\notin \bR$ and $\al,\be,\ga\in (0,1)$
such that $\ga\neq (1-\al)\be$ and $1-\ga\neq (1-\al)(1-\be)$. Then 
Corollary~\ref{cor-gs} implies that if 
$$\al \Phi(z_1)+(1-\al)\Phi(\be z_2+(1-\be)z_3)\le 
\ga \Phi(z_2)+(1-\ga)\Phi(z_3)$$ 
for all convex functions $\Phi:\bC\rightarrow \bR$ then $\ga>(1-\al)\be$, 
$1-\ga>(1-\al)(1-\be)$, and the inequality
$$\al \Phi(z_1)+(1-\al)(\be \Phi(z_2)+(1-\be)\Phi(z_3))\le \ga \Phi(z_2)+
(1-\ga)\Phi(z_3)$$
holds for any convex function $\Phi:\bC\rightarrow \bR$.
\end{remark}

\begin{remark}
Theorem~\ref{gen-S} generalizes several 
known results on vector and matrix majorization, many of which
seem in fact to be regularly rediscovered in various contexts,  
see e.g.~\cite{Da} and references therein. Two new interesting 
extensions of classical and multivariate majorization have 
recently been introduced and studied in \cite{Ma2}.
\end{remark}

\section{Geometry of equilibrium points: statement of main 
results}\label{s2} 

In this section we formulate our main results. These describe the distribution 
of equilibrium points of logarithmic potentials for discrete configurations 
of positive charges (Theorems~\ref{t-1}--\ref{t-4} and 
Theorems~\ref{t-7}--\ref{t-8}) and the Hausdorff geometry of the 
extended set of equilibrium points for finite charge configurations 
(Theorems~\ref{t-5}--\ref{t-6}).

\subsection{Weighted majorization relations for finite charge 
configurations}\label{s21}

Let $z_i$, $1\le i\le n$, be distinct points in the complex plane carrying 
positive charges $a_i$, $1\le i\le n$, respectively, where $n\ge 2$. 
Throughout 
this section we assume wlog that the $a_i$ are normalized so that
$\sum_{i=1}^{n}a_i=1$. The logarithmic potential generated by these charged 
particles is given by
\begin{equation}\label{pot}
U(z)=\sum_{i=1}^{n}a_{i}\log\left|1-\frac{z}{z_i}\right|
\end{equation}
and the resulting electrostatic force is
\begin{equation}\label{f}
f(z)=\sum_{i=1}^{n}\frac{a_i}{z-z_i}.
\end{equation}

The class of functions of the form~\eqref{f} is essentially 
the same as Sz-Nagy's class of generalized derivatives \cite{Sz}. 
As explained in \S \ref{s0}, the critical points of $U(z)$ coincide with 
the equilibrium points of the electrostatic field, that is, the zeros of 
$f(z)$. By an argument reminiscent of the Gauss-Lucas theorem one 
can see that all equilibrium points lie in $\co(z_1,\ldots,z_n)$, 
cf., e.g., Proposition 3.1 in \cite{Ma2}. Remarkable generalizations of this
simple fact were recently obtained in {\em loc.~cit.}, where it was shown 
that these points actually satisfy a whole sequence of majorization relations.
Theorems~\ref{t-1}--\ref{t-3} below further extend the aforementioned results.

{\allowdisplaybreaks
\begin{notation}
Let $w_1,\ldots,w_{n-1}$ be the critical points of $U(z)$ counted with 
multiplicity and recall the set $\calA_m$ defined in~\eqref{sp-n}. 
For $\la\in\bC\setminus\{0\}$ and $\mu\in\bC$ let
\begin{align}\label{not-1}
&W(\la,\mu)=(\la w_1+\mu,\ldots,\la w_{n-1}+\mu)\in\bC^{n-1},\\
&Z(\la,\mu)=(\la z_1+\mu,\ldots,\la z_n+\mu)\in\bC^{n},\notag\\
&\ba=\left(\frac{1}{n-1},\ldots,\frac{1}{n-1}\right)\in\calA_{n-1},\quad
\bb=\left(\frac{1-a_1}{n-1},\ldots,\frac{1-a_n}{n-1}\right)\in\calA_{n}.\notag
\end{align}
\end{notation}}

\begin{theorem}\label{t-1}
There exists a matrix $R\in \bM_{n-1,n}^{rs}$ such that 
$W(\la,\mu)^T=RZ(\la,\mu)^T$ and $\bb=\ba R$ for any 
$\la\in\bC\setminus\{0\}$ and $\mu\in\bC$, so that 
$(W(\la,\mu),\ba)\pr (Z(\la,\mu),\bb)$. Equivalently, the inequality
\begin{equation}\label{f-t-1}
\sum_{j=1}^{n-1}\Phi(\la w_j+\mu)\le \sum_{i=1}^{n}(1-a_i)\Phi(\la z_i+\mu)
\end{equation}
is satisfied by any (continuous) convex function $\Phi:\bC\to \bR$.
\end{theorem}

\begin{remark}\label{pere-mala}
The special case of~\eqref{f-t-1} when $\la=1$, $\mu=0$ and $a_i=\frac{1}{n}$,
$1\le i\le n$, is already a considerable improvement of the Gauss-Lucas 
theorem and was originally conjectured by de Bruijn and 
Springer \cite{BrS}. This conjecture was recently proved 
by S.~Malamud \cite{Ma1,Ma2} and R.~Pereira \cite{P} independently of each 
other.
\end{remark}

\begin{remark}\label{r-3}
Theorem~\ref{t-1} implies in particular that the following ordinary
multivariate majorization relation holds in $\bR^2$:
$$(w_1,\ldots,w_{n-1},\ze)\pr (z_1,\ldots,z_{n-1},z_n),$$
where $\ze=\sum_{i=1}^{n}a_{i}z_i$ is the $(a_1,\ldots,a_n)$-barycenter of 
$Z(1,0)$ (cf.~Remarks~\ref{bary} and~\ref{ordi}). This relation was previously
obtained in \cite[Proposition 4.3]{Ma2} and \cite[Theorem 5.4]{P}. Note though
that if $n\ge 3$ then even in the case of equal charges the above relation is 
in general weaker than the corresponding one given by Theorem~\ref{t-1} since 
the $z_i$ are linearly dependent (compare with Corollary~\ref{cor-gs}).
\end{remark}

As we shall now see, \eqref{f-t-1} is actually 
but one relation in a sequence of $\binom{n}{2}$
weighted multivariate majorization relations satisfied by the zeros and poles 
of $f(z)$. Moreover, the equilibrium points of $U(z)$ are in fact uniquely
determined by $n-1$ such relations. To formulate these results we need some 
new notation. 

\begin{notation}\label{new-not2}
Given $k\in \bN$ and $1\le m\le k$
let $\Pi_{k,m}$ denote the $m$th elementary symmetric function on $k$ symbols.
For $\la\in\bC\setminus \{0\}$, $\mu\in\bC$, $k\in\{1,\ldots,n-1\}$ and
$1\le m\le k$ set
{\allowdisplaybreaks
\begin{align}\label{not-4}
&W_{m}^{[k]}(\la,\mu)
=\Big(\om_{m}[r_1,\ldots,r_k](\la,\mu)\Big)_{1\le r_1<\ldots<r_k\le n-1}
\in \bC^{\binom{n-1}{k}},\\
&\text{where }\,\om_{m}[r_1,\ldots,r_k](\la,\mu)
=\Pi_{k,m}(\la w_{r_1}+\mu,\ldots,\la w_{r_k}+\mu),\notag\\
&Z_{m}^{[k]}(\la,\mu)
=\Big(\ze_{m}[s_1,\ldots,s_k](\la,\mu)\Big)_{1\le s_1<\ldots<s_k\le n}
\in \bC^{\binom{n}{k}},\notag\\
&\text{where }\,\ze_{m}[s_1,\ldots,s_k](\la,\mu)
=\Pi_{k,m}(\la z_{s_1}+\mu,\ldots,\la z_{s_k}+\mu),\notag\\
&\ba^{[k]}=
\left(\binom{n-1}{k}^{-1},\ldots,\binom{n-1}{k}^{-1}\right)
\in \calA_{\binom{n-1}{k}},\notag\\
&\bb^{[k]}=\left(\binom{n-1}{k}^{-1}
\!\!\left(1-\sum_{i=1}^{k}a_{s_i}\right)\right)_{1\le s_1<\ldots<s_k\le n}
\in \calA_{\binom{n}{k}}.\notag
\end{align}
Below we shall always assume that the components of $W_{m}^{[k]}(\la,\mu)$, 
$Z_{m}^{[k]}(\la,\mu)$ and $\bb^{[k]}$ are arranged 
increasingly with respect to the lexicographic order on the set of all indices 
consisting of $k$-tuples of increasing numbers in $\{1,\ldots,n\}$. 
In particular,} for any $1\le m\le \binom{n}{k}$ the $m$th coordinates
of $Z_{m}^{[k]}(\la,\mu)$ and $\bb^{[k]}$ are indexed by the same $k$-tuple 
$(s_1,\ldots,s_k)$ satisfying $1\le s_1<\ldots<s_k\le n$. Note also that with 
this ordering~\eqref{not-1} and~\eqref{not-4} imply that $\ba^{[1]}=\ba$,
$\bb^{[1]}=\bb$, $W_{1}^{[1]}(\la,\mu)=W(\la,\mu)$ and 
$Z_{1}^{[1]}(\la,\mu)=Z(\la,\mu)$. 
\end{notation}

Theorem 3.11 in \cite{Ma2} describes all possible sets of equilibrium points 
for logarithmic potentials generated by a fixed configuration of particles 
that are allowed to carry arbitrary positive charges. A natural related 
inverse problem 
is to describe geometrically the equilibrium points of a logarithmic potential
associated to a given set of particles with prescribed positive charges. 
The following theorem solves this problem and provides a characterization of 
the zeros of the 
logarithmic potential $U(z)$ defined in~\eqref{pot} by means of 
$n-1$ weighted majorization relations. 

\begin{theorem}\label{t-3}
The following conditions are equivalent: 
\begin{itemize}
\item[(i)] The zeros of $f(z)$ counted with multiplicity are 
$w_1,\ldots,w_{n-1}$, that is, 
\begin{equation}\label{cont}
\sum_{i=1}^{n}a_{i}\prod_{\substack{j=1\\j\neq i}}^{n}(z-z_j)
=\prod_{k=1}^{n-1}(z-w_k),\quad z\in\bC.
\end{equation}
\item[(ii)] There exist $\la\in\bC\setminus\{0\}$ and $\mu\in \bC$ such that
$$\left(W_{k}^{[k]}(\la,\mu),\ba^{[k]}\right)\pr 
\left(Z_{k}^{[k]}(\la,\mu),\bb^{[k]}\right),\quad 1\le k\le n-1.$$
\item[(iii)] There exist $\la\in\bC\setminus\{0\}$ and $\mu\in \bC$ such that
for any $k\in\{1,\ldots,n-1\}$ one can find an 
$\binom{n-1}{k}\times \binom{n}{k}$ row stochastic matrix $R_{k}(\la,\mu)$ 
satisfying
$$W_{k}^{[k]}(\la,\mu)^T=R_{k}(\la,\mu)Z_{k}^{[k]}(\la,\mu)^T\text{ and }
\bb^{[k]}=\ba^{[k]}R_{k}(\la,\mu).$$
\item[(iv)] There exist $\la\in\bC\setminus\{0\}$ and $\mu\in \bC$ such that
for any $k\in\{1,\ldots,n-1\}$ the inequality 
\begin{equation*}
\begin{split}
\sum_{1\le r_1<\ldots<r_k\le n-1}&\Phi
\Big(\om_{k}[r_1,\ldots,r_k](\la,\mu)\Big)\\
\le
&\sum_{1\le s_1<\ldots<s_k\le n}\left(1-\sum_{i=1}^{k}a_{s_i}\right)
\Phi\Big(\ze_{k}[s_1,\ldots,s_k](\la,\mu)\Big)
\end{split}
\end{equation*}
holds for all (continuous) convex functions $\Phi:\bC\to \bR$.
\item[(v)] Either of conditions (ii)--(iv) is true for $\la=1$ and $\mu=0$.
\end{itemize}
If the above conditions are fulfilled then for any 
$k\in\{1,\ldots,n-1\}$ there exists an 
$\binom{n-1}{k}\times \binom{n}{k}$ row stochastic matrix $R_{k}$ such that
\begin{equation}\label{inv-m}
W_{k}^{[k]}(\la,\mu)^T=R_{k}Z_{k}^{[k]}(\la,\mu)^T\text{ and }
\bb^{[k]}=\ba^{[k]}R_{k}
\end{equation}
for all $\la\in\bC\setminus\{0\}$ and $\mu\in \bC$.
\end{theorem}

\begin{remark}\label{dBS}
Theorem~\ref{t-3} extends \cite[Theorem 4.7]{Ma2} to the case of arbitrary 
(not necessarily equal) charges. If $a_i=\frac{1}{n}$, $1\le i\le n$, then 
Theorem~\ref{t-3} reduces to the characterization of the zeros
and critical points of degree $n$ complex polynomials by means of $n-1$ 
majorization relations which was given in \cite[Theorem 4.10]{Ma2}. 
\end{remark}

From Theorem~\ref{t-3} we deduce a number of corollaries.

\begin{corollary}\label{t-4}
For any $k\in\{1,\ldots,n-1\}$ there exists an 
$\binom{n-1}{k}\times \binom{n}{k}$ row stochastic matrix $R_{k}$ such that
\begin{equation}\label{gen-m}
W_{m}^{[k]}(\la,\mu)^T=R_{k}Z_{m}^{[k]}(\la,\mu)^T\text{ and }
\bb^{[k]}=\ba^{[k]}R_{k}
\end{equation}
whenever $1\le m\le k$, $\la\in\bC\setminus\{0\}$ and $\mu\in \bC$. In 
particular, for all such parameters $k$, $m$, $\la$, $\mu$ and any convex 
function $\Phi:\bC\to \bR$ the following inequality holds: 
\begin{equation}\label{g-conv-f}
\begin{split}
\sum_{1\le r_1<\ldots<r_k\le n-1}&\Phi
\Big(\om_{m}[r_1,\ldots,r_k](\la,\mu)\Big)\\
\le
&\sum_{1\le s_1<\ldots<s_k\le n}\left(1-\sum_{i=1}^{k}a_{s_i}\right)
\Phi\Big(\ze_{m}[s_1,\ldots,s_k](\la,\mu)\Big).
\end{split}
\end{equation}
\end{corollary}

\begin{remark}
Note that~\eqref{gen-m} is {\em \`a priori} 
stronger than~\eqref{g-conv-f} since it 
asserts that there exist matrices $R_k$ of the desired type that 
simultaneously realize the weighted majorization relations between the 
pairs $\left(W_{m}^{[k]}(\la,\mu),\ba^{[k]}\right)$ and 
$\left(Z_{m}^{[k]}(\la,\mu),\bb^{[k]}\right)$ for {\em all} 
admissible parameters $m$, $\la$ and $\mu$. 
\end{remark}

Setting $m=1$ in Corollary~\ref{t-4} yields the following generalization of
Theorem~\ref{t-1}, which corresponds to the case when $k=1$ in~\eqref{f-t-2}
below.

\begin{corollary}\label{t-2}
For any $k\in\{1,\ldots,n-1\}$ there exists an 
$\binom{n-1}{k}\times \binom{n}{k}$ row stochastic matrix
$R_k$ such that 
$$W_{1}^{[k]}(\la,\mu)^T=R_{k}Z_{1}^{[k]}(\la,\mu)^T\text{ and }
\bb^{[k]}=\ba^{[k]}R_{k}$$ 
for any $\la\in\bC\setminus\{0\}$ and $\mu\in \bC$. Thus
$$\left(W_{1}^{[k]}(\la,\mu),\ba^{[k]}\right)\pr 
\left(Z_{1}^{[k]}(\la,\mu),\bb^{[k]}\right),\quad 1\le k\le n-1,$$ 
or equivalently, for all convex functions $\Phi:\bC\to \bR$ and 
$k\in\{1,\ldots,n-1\}$ one has
\begin{equation}\label{f-t-2}
\begin{split}
\sum_{1\le r_1<\ldots<r_k\le n-1}&\Phi
\Big(\om_{1}[r_1,\ldots,r_k](\la,\mu)\Big)\\
&\le
\sum_{1\le s_1<\ldots<s_k\le n}\left(1-\sum_{i=1}^{k}a_{s_i}\right)
\Phi\Big(\ze_{1}[s_1,\ldots,s_k](\la,\mu)\Big).
\end{split}
\end{equation}
\end{corollary}

\begin{remark}\label{la0}
The condition $\la\neq 0$ may actually be omitted from the above 
inequa\-lities. However, if $\la=0$ the vectors defined in~\eqref{not-4} 
are independent of the $z_i$ and 
$w_j$ and the aforementioned inequalities become 
equalities. Thus the corresponding weighted majorization relations 
are trivially satisfied in this case. 
\end{remark}

Taking $\la=1$, $\mu=0$ and $\Phi(z)=|z|^{\al}$ with $\al\ge 1$ in 
Corollary~\ref{t-4} we obtain the following moment inequalities for the 
charged particles and the equilibrium points.

\begin{corollary}\label{cor-mom}
For any $k\in\{1,\ldots,n-1\}$, $1\le m\le k$, and $\al\ge 1$ one has
\begin{equation*}
\begin{split}
\sum_{1\le r_1<\ldots<r_k\le n-1}&
\left|\Pi_{k,m}(w_{r_1},\ldots,w_{r_k})\right|^{\al}\\
&\le\sum_{1\le s_1<\ldots<s_k\le n}
\left|\Pi_{k,m}(z_{s_1},\ldots,z_{s_k})\right|^{\al},
\end{split}
\end{equation*}
where $\Pi_{k,m}$ denotes as before the $m$th elementary symmetric function on
$k$ symbols.
\end{corollary}

In particular, by setting $a_i=\frac{1}{n}$, $1\le i\le n$, and $m=1$ in 
Corollary~\ref{cor-mom} we get 
\begin{equation}\label{r-s}
\binom{n-1}{k}^{-1}\sum_{1\le r_1<\ldots<r_k\le n-1}
\left|\sum_{i=1}^{k}w_{r_i}\right|^{\al}
\le
\binom{n}{k}^{-1}\sum_{1\le s_1<\ldots<s_k\le n}
\left|\sum_{i=1}^{k}z_{s_i}\right|^{\al}
\end{equation}
for all $\al\ge 1$ and $k\in\{1,\ldots,n-1\}$. Note that for $k=2$ and 
$\al\ge 1$ the inequalities in~\eqref{r-s} are quite similar to the 
inequalities 
$$\binom{n-1}{2}^{-1}\sum_{1\le r_1<r_2\le n-1}
\left|w_{r_1}-w_{r_2}\right|^{\al}
\le \binom{n}{2}^{-1}\sum_{1\le s_1<s_2\le n}
\left|z_{s_1}-z_{s_2}\right|^{\al}$$
that were conjectured by Rahman and Schmeisser in \cite[Remark 2.3.9]{RS}. In 
\S \ref{s-mom} we conjecture a whole hierarchy of weighted 
de Bruijn-Springer relations and inertia laws, see Conjecture~\ref{gen-mom}
below. The latter actually contains as special cases all types of moment 
inequalities like those that we just described. 

\subsection{Relative Hausdorff geometry}\label{s212}

The metric function between (closed) subsets of Euclidean space introduced by 
Hausdorff makes use of the following notions.

\begin{definition}\label{d-HT}
The {\em relative} (alternatively, {\em directed} or {\em oriented}) 
{\em Hausdorff distance} from a subset $\Om_1$ of $\bC$ to another 
such subset $\Om_2$ is given by
$$h(\Om_1,\Om_2)=\sup_{w\in\Om_1}\text{dist}(w,\Om_2)
=\sup_{w\in\Om_1}\inf_{z\in\Om_2}|w-z|.$$
Other frequent names and notations are the {\em forward Hausdorff distance} 
from $\Om_1$ to $\Om_2$, $h_f(\Om_1,\Om_2):=h(\Om_1,\Om_2)$, and the 
{\em backward Hausdorff distance} from $\Om_1$ to 
$\Om_2$, $h_b(\Om_1,\Om_2):=h(\Om_2,\Om_1)$. The 
{\em (symmetrized) Hausdorff distance} between $\Om_1$ and $\Om_2$ is then 
given by
$$H(\Om_1,\Om_2)=\max\big(h(\Om_1,\Om_2),h(\Om_2,\Om_1)\big).$$ 
\end{definition}

Let $\calZ=\{z_1,\ldots,z_n\}$ be a set of charged particles and
$\calW=\{w_1,\ldots,w_{n-1}\}$ be the corresponding set of equilibrium points 
of the logarithmic potential defined in~\eqref{pot}.
Theorems~\ref{t-1}--\ref{t-3} yield a 
series of relations between the sets $\calW$ and $\calZ$ when all their 
elements are considered simultaneously. These relations go much beyond the 
obvious inclusion $\calW\subset \co(\calZ)$ and suggest that 
no equilibrium point can actually lie too far away from the 
particles. Our next result gives a precise meaning to this intuitive 
fact and 
shows that the relative Hausdorff distance $h(\calW,\calZ)$ is 
dominated by the {\em $\ba$-weighted standard deviation} of the 
particles. The latter is given by the square root of the 
{\em $\ba$-weighted variance} of $\bz=(z_1,\ldots,z_n)$, which we define as 
\begin{equation}\label{wdev}
\si_2(\bz;\ba)^2=\inf_{\al\in\bC}\sum_{i=1}^{n}a_i|z_i-\al|^2
=\sum_{i=1}^{n}a_{i}|z_i-\ze|^2,
\end{equation}
where $\ba=(a_1,\ldots,a_n)\in\calA_n$ and $\ze$ stands as before for the 
$\ba$-weighted barycenter 
$\sum_{i=1}^{n}a_{i}z_i$ (cf.~Remark~\ref{r-3}). We shall actually prove
a stronger statement involving the {\em extended set of equilibrium points} 
of the given system of positive charges, i.e., the set 
\begin{equation}\label{ext}
\calW_e=\calW\cup \{\ze\}=\{w_1,\ldots,w_{n-1},\ze\}.
\end{equation}

\begin{theorem}\label{t-5}
In the above notation one has 
$h(\calW,\calZ)\leq h(\calW_e,\calZ)\leq \si_2(\bz;\ba)$.
\end{theorem}

\begin{remark}\label{optim}
The bound given in Theorem~\ref{t-5} is sharp, as one can see for instance 
by placing $n$ equal charges at the vertices of a regular $n$-gon. 
\end{remark}

A natural question in this context is whether there is any analog of 
Theorem~\ref{t-5} for the relative Hausdorff distance from $\calZ$ 
to $\calW_e$. The following
theorem provides an affirmative answer to this question when all 
particles are assumed to lie on the same line and 
shows that in this case the Hausdorff geometry of the set of 
charged particles and 
the extended set of equilibrium points is controlled by the weighted standard 
deviation of the particles.

\begin{theorem}\label{t-6}
If $z_1,\ldots,z_n$ are collinear then 
$H(\calZ,\calW_e)\leq \si_2(\bz;\ba)$.
\end{theorem}  

\begin{remark}\label{optim-e}
Note that Theorem~\ref{t-6} is sharp and that a similar result cannot hold if
$\calW_e$ is replaced by $\calW$. To see this one may simply consider two 
distinct points in the plane carrying unequal positive charges.
\end{remark}

\subsection{Infinite discrete distributions of charge}\label{s22}

As we point out in \S\ref{s42}, several difficulties occur when trying to 
extend the results of \S \ref{s21} to logarithmic potentials generated by 
an infinite discrete set of positive charges. To begin with, in this case 
already the fundamental 
question dealing with the {\em existence} of equilibrium points is quite 
delicate and, as a matter of fact, still an open problem 
(cf.~\cite{CER,ELR,LR}). 

Nevertheless, it turns out that in many situations one 
can obtain new information on the geometry of the electrostatic field through 
natural infinite-dimensional versions of Theorem~\ref{t-1}. As we shall now 
explain, this is for instance the case for collinear positive charges and, 
more generally, for essentially any bounded discrete set of positively 
charged particles with finite total charge for which equilibrium points are 
known to exist. Let 
$\{a_i\}_{i\in\bN}$ be a sequence of positive numbers and let 
$\{z_i\}_{i\in\bN}$ be a sequence of distinct complex numbers satisfying the 
conditions
\begin{equation}\label{i-hyp}
a_i>0\text{ and }|z_i|<\rho\text{ for }i\in\bN,\quad \sum_{i=1}^{\infty}a_i=1,
\quad \lim_{i\to\infty}z_i=\rho,
\end{equation}
where $\rho$ is a fixed positive number. Using the same electrostatic 
interpretation as before -- that is, placing
a charge $a_i$ at $z_i$ for $i\in \bN$ -- we get the resulting force  
\begin{equation}\label{i-f}
f(z)=\sum_{i=1}^{\infty}\frac{a_i}{z-z_i}
\end{equation}
and an associated logarithmic potential given by
\begin{equation*}
U(z)=\sum_{i=1}^{\infty}a_{i}\log\left|1-\frac{z}{z_i}\right|.
\end{equation*}
Set $\bD(\rho)=\{z\in\bC: |z|<\rho\}$ and 
$\bD_{f}(\rho)=\bD(\rho)\setminus \{z_i\}_{i\in\bN}$. By~\eqref{i-hyp} one has 
\begin{equation*}
\sum_{i=1}^{\infty}\frac{a_i}{|z_i|}<\infty,
\end{equation*}
so that the series in the right-hand side 
of~\eqref{i-f} converges absolutely for all $z\in \bD_{f}(\rho)$ and thus 
$f$ is 
a meromorphic function in $\bD(\rho)$. One can see in similar fashion 
that $U$ 
is subharmonic in $\bD(\rho)$. We have the following analog of 
Theorem~\ref{t-1}.

\begin{theorem}\label{t-7}
If $\{z_i\}_{i\in\bN}$ is a real sequence satisfying~\eqref{i-hyp} then the 
function $f$ defined by~\eqref{i-f} has an infinite discrete set of real 
zeros $\{w_j\}_{j\in\bN}$ and the sequences $\{z_i\}_{i\in\bN}$ and 
$\{w_j\}_{j\in\bN}$ interlace on the real axis. Moreover, there exists an 
infinite matrix 
$S=(s_{ij})_{i,j=1}^{\infty}$ such that for all $(i,j)\in\bN^2$ one has
\begin{equation*}
s_{ij}\ge 0,\quad \sum_{k=1}^{\infty}s_{ik}=1,\quad 
\sum_{l=1}^{\infty}s_{lj}=1-a_j,\quad w_i=\sum_{k=1}^{\infty}s_{ik}z_k.
\end{equation*}
In particular, for any $\la\in\bC\setminus\{0\}$ and $\mu\in\bC$ the inequality
\begin{equation*}
\sum_{j=1}^{\infty}\Phi(\la w_j+\mu)\le 
\sum_{i=1}^{\infty}(1-a_i)\Phi(\la z_i+\mu)
\end{equation*}
holds for any convex function $\Phi:\bC\to\bR$.
\end{theorem} 

Our last main result is a natural analog of Theorem~\ref{t-7} for 
all bounded discrete configurations of positive charges for 
which equilibrium points are known to exist.

\begin{theorem}\label{t-8}
Let $\{z_i\}_{i\in\bN}$ be a sequence of distinct complex numbers 
satisfying~\eqref{i-hyp} and let $f$ be the meromorphic function in 
$\bD(\rho)$ 
defined by~\eqref{i-f}. If $f$ has a non-empty (discrete) set of 
zeros $\{w_j\}_{j\in\bN}$ then the inequality
$$ \sum_{j=1}^{\infty}\Phi(\la w_j+\mu)\le 
\sum_{i=1}^{\infty}(1-a_i)\Phi(\la z_i+\mu)$$
holds for all $\la\in\bC\setminus\{0\}$, $\mu\in\bC$ and any nonnegative 
convex function $\Phi:\bC\to\bR$.
\end{theorem} 

\begin{remark}\label{expl}
The fact that inequality~\eqref{bes-2} in \S \ref{s33} may be strict 
explains why the 
inequality in Theorem~\ref{t-8} holds in general only for {\em nonnegative} 
convex functions.
The connection between majorization and such inequalities for nonnegative 
convex functions was first found in \cite{FH}.
\end{remark}

\begin{remark}\label{kat}
Theorem~\ref{t-8} is reminiscent of Kato's bound for the
variation of discrete spectra of Hilbert space selfadjoint operators under 
compact perturbations \cite{Ka1}.
\end{remark}

\section{Proof of main results}\label{s3}

\subsection{Proof of Theorems~\ref{t-1}--\ref{t-3}}\label{s31}

We shall use various methods from multilinear algebra and matrix analysis. 
These require a few notations and preliminary results and we refer 
to~\cite{Bh} and~\cite{MO} for background material. 

\begin{notation}\label{not-proofs}
Given $m\in\bN$ 
and $1\le k\le m$
let $$Q_{k,m}=\{\bi=(i_1,\ldots,i_k)\mid 1\le i_1<\ldots<i_k\le m\}.$$
If $B=(b_{ij})$ is an $m\times m$ matrix and $\bi,\bj\in Q_{k,m}$
we denote by
$$B(\bi,\bj)=B\binom{i_1,\ldots,i_k}{j_1,\ldots,j_k}$$
the $k\times k$ submatrix of $B$ lying in rows $i_1,\ldots,i_k$ and columns
$j_1,\ldots,j_k$. Below we shall always assume that the elements of $Q_{k,m}$
are arranged in lexicographic order, i.e., if $\bi$ and $\bj$ are distinct
$k$-tuples in $Q_{k,m}$ then $\bi\ge \bj$ if the first non-zero term in the
sequence $i_1-j_1,\ldots,i_k-j_k$ is positive.
\end{notation}

We fix an $n$-dimensional 
complex Hilbert space $\calH$ with unitarily invariant scalar 
product $\langle\cdot,\!\cdot\rangle$ and identity operator 
$I=I_{\calH}=I_n\in L(\calH)$, where $L(\calH)$ is the
set of all linear operators on $\calH$. Let $A\in L(\calH)$ be a normal
operator with spectrum $\Si(A)=\{z_1,\ldots,z_n\}$ and choose an 
orthonormal basis $(\mathbf{e}_1,\ldots,\mathbf{e}_n)$ of $\calH$ such that 
$A\bbe_i=z_i\bbe_i$ for $1\le i\le n$. Define a unit vector
\begin{equation}\label{v-n}
\mathbf{v}_{n}=\sum_{i=1}^{n}\sqrt{a_i}\bbe_i
\end{equation}
and let $P$ be the orthogonal projection on the subspace 
$\calK:=\mathbf{v}_{n}^{\perp}$ of $\calH$. The operator
$$A'=PAP|_{\calK}\in L(\calK)$$
is called the $P$-{\em compression} of $A$ \cite{D}. Recall the function $f$ 
defined in~\eqref{f}. The following statement may be found in 
e.g.~\cite[formula (3.3)]{Ma2} and \cite[Lemma 2.2]{P}. 

\begin{lemma}\label{g-comp}
With the above notations one has
$$\left\langle (A-zI_{\calH})^{-1}\mathbf{v}_n,\mathbf{v}_n\right\rangle=
\frac{\det(A'-zI_{\calK})}{\det(A-zI_{\calH})}=-f(z)\text{ for }z\in\bC
\setminus \Si(A).$$
\end{lemma} 

\begin{proof}
Let $(\mathbf{v}_1,\ldots,\mathbf{v}_{n-1})$ be an orthonormal basis of 
$\calK$. The
matrix representation of $A'$ in the basis $(\mathbf{v}_1,\ldots,
\mathbf{v}_{n-1})$ is given by 
the $(n-1)\times (n-1)$ upper left-hand principal submatrix of the matrix 
representation of $A$ in the orthonormal basis $(\mathbf{v}_1,\ldots,
\mathbf{v}_{n-1},\mathbf{v}_n)$ of 
$\calH$. For any $z\in \bC
\setminus \Si(A)$ the $(n,n)$ entry of the 
matrix representation of the (normal) operator $(A-zI_{\calH})^{-1}$ in 
the basis $(\mathbf{v}_1,\ldots,\mathbf{v}_{n-1},\mathbf{v}_n)$ is given on 
the one hand by
\begin{equation*}
\left\langle (A-zI_{\calH})^{-1}\mathbf{v}_n,\mathbf{v}_n\right\rangle 
=\left\langle \sum_{i=1}^{n}\frac{\langle \mathbf{v}_n,\mathbf{e}_i\rangle}
{z_i-z}\mathbf{e}_i,
\sum_{i=1}^{n}\langle \mathbf{v}_n,\mathbf{e}_i\rangle \mathbf{e}_i
\right\rangle=-f(z).
\end{equation*}
On the other hand, by Cramer's rule the $(n,n)$ entry of the matrix 
representation of $(A-zI_{\calH})^{-1}$ in 
the basis $(\mathbf{v}_1,\ldots,\mathbf{v}_{n-1},\mathbf{v}_n)$ of 
$\calH$ is given by the cofactor of the $(n,n)$ 
entry of the matrix representation of $(A-zI_{\calH})$ in the same 
basis. Thus 
$$\left\langle (A-zI_{\calH})^{-1}\mathbf{v}_n,\mathbf{v}_n\right\rangle=
\frac{\det(A'-zI_{\calK})}{\det(A-zI_{\calH})}$$
for $z\in \bC
\setminus \Si(A)$, which proves the lemma.
\end{proof}

\begin{proof}[Proof of Theorem~\ref{t-1}]
By Lemma~\ref{g-comp} the zeros $w_1,\ldots,w_{n-1}$ of $f$ are the same
as the eigenvalues of the $P$-compression $A'$. 
Let $(\bv_1,\ldots,\bv_{n-1})$ be an orthonormal 
basis of $\calK$ that triangularizes $A'$. Then 
\begin{equation}\label{f-1}
w_j=\bv_{j}^{*}A'\bv_j=\langle A'\bv_j,\bv_j\rangle=\langle A\bv_j,\bv_j\rangle
=\sum_{i=1}^{n}z_i|\langle \bv_j,\bbe_i\rangle|^2
\end{equation}
since $\bv_j=\sum_{i=1}^{n}\langle \bv_j,\bbe_i\rangle\bbe_i$ for $1\le j\le 
n-1$. Define the $(n-1)\times n$ matrix $R=(r_{ij})$ by setting $r_{ij}=
|\langle \bv_i,\bbe_j\rangle|^2$ for $1\le i\le n-1$ and $1\le j\le n$. Using 
the fact that $(\bv_1,\ldots,\bv_{n-1},\bv_n)$ is an orthonormal basis of 
$\calH$ one can easily check that
\begin{equation*}
\begin{split}
&\sum_{j=1}^{n}r_{ij}=\sum_{j=1}^{n}|\langle \bv_i,\bbe_j\rangle|^2=||\bv_i||^2
=1,\quad 1\le i\le n-1,\\
&\sum_{i=1}^{n-1}r_{ij}=\sum_{i=1}^{n-1}|\langle \bv_i,\bbe_j\rangle|^2=
||\bbe_j||^2-|\langle \bbe_j,\bv_n\rangle|^2=1-a_{j},\quad 1\le j\le n,
\end{split}
\end{equation*}
so that $R\in\bM_{n-1,n}^{rs}$ and $\bb=\ba R$. Note that by~\eqref{f-1} 
one has 
$$(w_1,\ldots,w_{n-1})^T=R(z_1,\ldots,z_n)^T.$$ 
Hence $(W(1,0),\ba)\pr (Z(1,0),\bb)$, which completes the proof
of the theorem in the case $\la=1$ and $\mu=0$. The general case follows 
from this one by composing $f$ with a non-singular affine transformation of 
the plane.
\end{proof}

\begin{remark}\label{normcompr}
The compression of a normal operator to a hyperplane 
is not necessarily normal. For instance, Fan and Pall \cite{FP} showed that 
principal submatrices of a given irreducible normal 
matrix $B$ are normal if and only if $B$ has collinear eigenvalues. 
\end{remark}  

\begin{remark}
The proof of Theorem~\ref{t-1} given above is a slight modification of the 
corresponding arguments in \cite{Ma2} and \cite{P}.
\end{remark} 

We next prove Theorem~\ref{t-3}. Clearly, conditions (ii), (iii) and
(iv) in this theorem are equivalent by Definition~\ref{def-w} and
Theorem~\ref{gen-S}. Let us first show that (v) $\Rightarrow$ (i). Set
$$q(z)=\prod_{j=1}^{n-1}(z-w_j),\quad p_{i}(z)
=\prod_{\substack{j=1\\j\neq i}}^{n}(z-z_j),\quad 1\le i\le n,$$
and note that condition (i) is equivalent to
\begin{equation*}
q^{(n-1-k)}(0)=\sum_{i=1}^{n}a_{i}p_{i}^{(n-1-k)}(0),\quad 1\le k\le n-1.
\end{equation*}
Elementary computations show that
\begin{equation*}
\begin{split}
&q^{(n-1-k)}(0)=(-1)^{k}(n-1-k)!
\sum_{1\le r_1<\ldots<r_k\le n-1}\prod_{j=1}^{k}w_{r_j},\\
&\sum_{i=1}^{n}a_{i}p_{i}^{(n-1-k)}(0)=(-1)^{k}(n-1-k)!
\sum_{1\le s_1<\ldots<s_k\le n}\left(1-\sum_{i=1}^{k}a_{s_i}\right)
\prod_{i=1}^{k}z_{s_i},
\end{split}
\end{equation*}
for $1\le k\le n-1$. Thus condition (i) amounts to saying that if 
$1\le k\le n-1$ then
\begin{equation}\label{id-k}
\sum_{1\le r_1<\ldots<r_k\le n-1}\prod_{j=1}^{k}w_{r_j}=
\sum_{1\le s_1<\ldots<s_k\le n}\left(1-\sum_{i=1}^{k}a_{s_i}\right)
\prod_{i=1}^{k}z_{s_i}.
\end{equation}
The fact that (v) $\Rightarrow$ (i) is now a 
consequence of the following lemma.

\begin{lemma}\label{easy}
If condition (v) in Theorem~\ref{t-3} holds then~\eqref{id-k} is true for 
$1\le k\le n-1$.
\end{lemma}

\begin{proof}
Using e.g.~condition (iii) in Theorem~\ref{t-3} with $\la=1$ and $\mu=0$ 
one gets
$$\ba^{[k]}W_{k}^{[k]}(1,0)^T=\bb^{[k]}Z_{k}^{[k]}(1,0)^T,\quad 
1\le k\le n-1,$$
which is the same as the desired conclusion.
\end{proof}

We shall also need the following extension to normal matrices of a 
well-known result of Schur for Hermitian matrices (cf.~\cite{Bh,MO}), see 
e.g.~\cite[Proposition 3.9 (a)]{Ma2}.

\begin{lemma}\label{s-lem}
If $T=(t_{ij})$ is a normal $n\times n$ matrix with eigenvalues 
$\tau_1,\ldots,\tau_n$ then
$$(t_{11},\ldots,t_{nn})\pr (\tau_1,\ldots,\tau_n),$$
i.e., $(t_{11},\ldots,t_{nn})^T=S(\tau_1,\ldots,\tau_n)^T$ for some 
doubly stochastic $n\times n$ matrix $S$.
\end{lemma}

\begin{proof}
Let $U$ be a unitary matrix such that 
$T=U\text{diag}(\tau_1,\ldots,\tau_n)U^{*}$ and set $S=U\circ \bar{U}$, where
$\circ$ denotes the Hadamard-Schur (entrywise) product of matrices and 
$\bar{U}$ is 
the complex conjugate of $U$. Then 
$(t_{11},\ldots,t_{nn})^T=S(\tau_1,\ldots,\tau_n)^T$. It remains to note that
$S$ is a unitarily stochastic hence a doubly stochastic matrix.
\end{proof}

To complete the proof of Theorem~\ref{t-3} it will be convenient 
to work with certain suitably chosen matrix representations of the operators
$A$ and $A'$ defined above. As in the proof of Theorem~\ref{t-1} let us fix 
an orthonormal 
basis $(\bv_1,\ldots,\bv_{n-1})$ of $\calK$ in 
which the matrix representation of the $P$-compression $A'$ of $A$ is upper
triangular. For simplicity, denote the matrix representation of $A$ in the
basis $(\bv_1,\ldots,\bv_{n-1},\bv_n)$ again by $A$ and let $A_i$, 
$1\le i\le n$, be the degeneracy one principal submatrix of $A$ obtained by 
deleting its $i$-th row and $i$-th column. Then clearly $A_n=A'$ and so by
Lemma~\ref{g-comp} the spectrum of $A_n$ is 
$\Si(A_n)=\{w_1,\ldots,w_{n-1}\}$, where 
$\Si(A_n)$ is viewed as a multiset whose elements occur as many times as their 
algebraic multiplicities as eigenvalues of $A_n$. Thus $A$ is an $n\times n$ 
normal matrix with $\Si(A)=\{z_1,\ldots,z_n\}$ such that 
\begin{equation}\label{form}
A=\begin{pmatrix}
A_{n} & \ast\\
\ast & \ze
\end{pmatrix},
\end{equation}
where $A_{n}$ is upper triangular with 
$\Si(A_n)=\{w_1,\ldots,w_{n-1}\}$ and
$\ze=\langle A\bv_n,\bv_n\rangle=\sum_{i=1}^{n}a_{i}z_i$ is the 
$(a_1,\ldots,a_n)$-barycenter of $Z(1,0)$ (cf.~Remark~\ref{r-3}). Note that
since the $w_j$ appear on the diagonal of $A_n$ Lemma~\ref{s-lem} actually
provides a new proof of Theorem~\ref{t-1}. Let $U=(u_{ij})$ be
the (unitary) change of basis matrix from $(\bv_1,\ldots\bv_n)$ to 
$(\bbe_1,\ldots,\bbe_n)$ so that we may write
\begin{equation}\label{dec}
A=U\text{diag}(z_1,\ldots,z_n)U^{*}.
\end{equation}

\begin{proof}[Proof of Theorem~\ref{t-3}]
Let us first show that if~\eqref{cont} holds and $k\in\{1,\ldots,n-1\}$
then one can find an $\binom{n-1}{k}\times\binom{n}{k}$ row stochastic 
matrix $R_k$ such that
\begin{equation}\label{want}
W_{k}^{[k]}(1,0)^T=R_{k}Z_{k}^{[k]}(1,0)^T
\text{ and }\bb^{[k]}=\ba^{[k]}R_k.
\end{equation}
This is clear for $k=1$ by Theorem~\ref{t-1} and so we may assume that 
$k\ge 2$. Let
$$A^{(k)}=\left(\det A(\bi,\bj)\right)_{\bi.\bj\in Q_{k,n}}$$
be the $k$-th compound matrix of $A$ acting on the $k$-th Grassmann power
$\wedge^{k}\calH$ (recall that the elements of $Q_{k,n}$ 
and therefore also the
$\binom{n}{k}\times\binom{n}{k}$ determinants appearing in $A^{(k)}$ 
are arranged in lexicographic order). The Binet-Cauchy 
formula for compound matrices \cite[Theorem 19.F.2]{MO} and~\eqref{dec}
imply that
$$A^{(k)}=U^{(k)}\text{diag}(z_1,\ldots,z_n)^{(k)}U^{*(k)},$$
which by an argument similar to the proof of Lemma~\ref{s-lem} leads to
\begin{equation}\label{stoc}
\left(\det A(\bi,\bi)\right)_{\bi\in Q_{k,n}}^{T}=S_k Z_{k}^{[k]}(1,0)^T,
\quad S_k=U^{(k)}\circ \bar{U}^{(k)},
\end{equation}
where $\circ$ denotes as before the Hadamard-Schur product. Note that $S_k$ 
is an 
$\binom{n}{k}\times\binom{n}{k}$ doubly stochastic matrix since $U^{(k)}$ is 
unitary. Clearly, the diagonal elements of $A^{(k)}$ corresponding to the 
$k\times k$ 
principal minors of $A_n$ are of the form
$\det A(\bi,\bi)$, $\bi\in Q_{k,n-1}$. There are a total of $\binom{n-1}{k}$
such diagonal elements and these must actually coincide with the coordinates 
of the vector $W_{k}^{[k]}(1,0)^T$ since $A_n$ is upper triangular. Let $R_k$
be the $\binom{n-1}{k}\times\binom{n}{k}$ row stochastic 
matrix obtained from $S_k$ by deleting all elements $S_k(\bi,\bj)$ with 
$\bi\in Q_{k,n}\setminus Q_{k,n-1}$. We show that $R_k$ 
satisfies~\eqref{want}. For this one has to prove that for any 
$\bj=(j_1,\ldots,j_k)\in Q_{k,n}$ the identity
$$\sum_{1\le i_1<\ldots<i_{k-1}\le n-1}
S_{k}\binom{i_1,\ldots,i_{k-1},n}{j_1,\ldots,j_k}=\sum_{l=1}^{k}a_{j_l}$$
holds, which is equivalent to
\begin{equation}\label{want1}
\sum_{1\le i_1<\ldots<i_{k-1}\le n-1}
\left|\det U\binom{i_1,\ldots,i_{k-1},n}{j_1,\ldots,j_k}\right|^2
=\sum_{l=1}^{k}a_{j_l}.
\end{equation}
or
\begin{equation}\label{want2}
\sum_{\bi\in Q_{k,n-1}}\left|\det U(\bi,\bj)\right|^2
=1-\sum_{l=1}^{k}a_{j_l}.
\end{equation}
By~\eqref{dec} one has $u_{nj}=\sqrt{a_j}$, $1\le j\le n$, so that by expanding
each determinant with respect to the last row we get 
$$\det U\binom{i_1,\ldots,i_{k-1},n}{j_1,\ldots,j_k}=\sum_{m=1}^{k}(-1)^{k+m}
\sqrt{a_{j_m}}
\det U\binom{i_1,\ldots,i_{k-1}}{j_1,\ldots,\widehat{j_m},\ldots,j_k}$$
and thus
\begin{multline*}
\left|\det U\binom{i_1,\ldots,i_{k-1},n}{j_1,\ldots,j_k}\right|^2\\
=\sum_{r,s=1}^{k}(-1)^{r+s}\sqrt{a_{j_r}a_{j_s}}
\det U\binom{i_1,\ldots,i_{k-1}}{j_1,\ldots,\widehat{j_r},\ldots,j_k}
\det \bar{U}\binom{i_1,\ldots,i_{k-1}}{j_1,\ldots,\widehat{j_s},\ldots,j_k}
\end{multline*}
for $1\le i_1<\ldots<i_{k-1}\le n-1$. Therefore
\begin{equation*}
\sum_{1\le i_1<\ldots<i_{k-1}\le n-1}
\left|\det U\binom{i_1,\ldots,i_{k-1},n}{j_1,\ldots,j_k}\right|^2\\
=\sum_{r,s=1}^{k}(-1)^{r+s}\sqrt{a_{j_r}a_{j_s}}\al_{r,s}(\bj),
\end{equation*}
where
\begin{equation*}
\al_{r,s}(\bj)=\sum_{1\le i_1<\ldots<i_{k-1}\le n-1}
\det U\binom{i_1,\ldots,i_{k-1}}{j_1,\ldots,\widehat{j_r},\ldots,j_k}
\det \bar{U}\binom{i_1,\ldots,i_{k-1}}{j_1,\ldots,\widehat{j_s},\ldots,j_k}.
\end{equation*}
Let us consider the $k\times (n-1)$ matrix
$$M=U\binom{1,\ldots,n-1}{j_1,\ldots,j_k}^T$$
and the $k\times k$ matrix $B=(b_{ij})=MM^{*}$ as well as the 
$k\times k$ matrix given by its $k-1$ 
Grassmann power $B^{(k-1)}=\left(b_{ij}^{(k-1)}\right)$. The Binet-Cauchy 
formula again implies that
$$B^{(k-1)}=M^{(k-1)}M^{*(k-1)}$$
from which we deduce that 
$$b_{k+1-r,k+1-s}^{(k-1)}=\al_{r,s}(\bj),\quad 1\le r,s\le k.$$
Thus the left-hand side of~\eqref{want1} equals
$$\sum_{r,s=1}^{k}(-1)^{r+s}\sqrt{a_{j_r}a_{j_s}}b_{k+1-r,k+1-s}^{(k-1)}$$
and by further manipulating this expression we obtain
$$\sum_{1\le i_1<\ldots<i_{k-1}\le n-1}
\left|\det U\binom{i_1,\ldots,i_{k-1},n}{j_1,\ldots,j_k}\right|^2
=\sum_{\bi\in Q_{k-1,k}}\det B(\bi,\bi)-k\det B$$
hence
$$\sum_{\bi\in Q_{k,n-1}}\left|\det U(\bi,\bj)\right|^2
=1-\sum_{\bi\in Q_{k-1,k}}\det B(\bi,\bi)+k\det B.$$
To show that the right-hand side of the above identity equals 
the right-hand side of~\eqref{want2} we proceed as follows. Define the vector 
$\bu=(\sqrt{a_{j_1}},\ldots,\sqrt{a_{j_k}})$. Since $U$ is 
unitary an easy computation shows that $B=I_k-\bu^T \bu$, so that $B$ depends
only on the numbers $a_{j_l}$, $1\le l\le k$. It follows that the expression
$$\sum_{1\le i_1<\ldots<i_k\le n-1}
\left|\det U\binom{i_1,\ldots,i_k}{j_1,\ldots,j_k}\right|^2
=\sum_{\bi\in Q_{k,n-1}}\left|\det U(\bi,\bj)\right|^2$$
also depends exclusively on the numbers $a_{j_l}$, $1\le l\le k$, and is 
therefore independent of the $z_i$, $1\le i\le n$. On the other hand, by
summing all coordinates corresponding to $\bi\in Q_{k,n-1}$ in each side
of~\eqref{stoc} we get
\begin{equation}\label{n-id-k}
\begin{split}
\sum_{1\le r_1<\ldots<r_k\le n-1}\prod_{j=1}^{k}&w_{r_j}\\
=\sum_{1\le s_1<\ldots<s_k\le n}&\left(\sum_{1\le i_1<\ldots<i_k\le n-1}
\left|\det U\binom{i_1,\ldots,i_k}{s_1,\ldots,s_k}\right|^2\right)
\prod_{i=1}^{k}z_{s_l}.
\end{split}
\end{equation}
Obviously, \eqref{cont} defines the $w_j$, $1\le j\le n-1$, as continuous
functions of the $z_i$, $1\le i\le n$. By letting $z_i\to 0$ for
$i\in \{1,\ldots,n\}\setminus \{j_1,\ldots,j_k\}$ and $z_{j_l}\to 1$ for
$1\le l\le k$ in both~\eqref{n-id-k} 
and~\eqref{id-k} and comparing coefficients in the resulting identities
we get~\eqref{want2}, which proves~\eqref{want}.

Finally, to prove the last statement of the theorem it is enough to consider 
the function $f(\la^{-1}(z-\mu))$, where $\la\in\bC\setminus\{0\}$
and $\mu\in\bC$. The zeros and 
poles of this function are the eigenvalues of 
$\la A_n+\mu$ and $\la A+\mu$, respectively. These matrices satisfy 
relations similar to~\eqref{form} and~\eqref{dec} while the unitary matrix 
$U$ is unchanged. 
One may therefore use the same arguments as above to get~\eqref{inv-m}. 
\end{proof}

\begin{remark}\label{cred-Ma}
Modulo some technical modifications that have to deal with unequal weights, 
the arguments given in the proof of Theorem~\ref{t-3} are the same as those
given by S.~Malamud in his proof of \cite[Theorem 4.7]{Ma2}.
\end{remark}

\begin{remark}\label{m-mat}
An interesting byproduct of the proof of Theorem~\ref{t-3} 
is the following identity for a certain type of $M$-matrices. Let 
$m\in \bN$, $\bx=(x_1,\ldots,x_m)\in (0,1)^m$ with 
$||\bx||^2=\sum_{i=1}^{m}x_{i}^2<1$ and define the $m\times m$ matrix
$X=I_m-\bx^T\bx$. Clearly, $\bx^T\bx$ is a rank one matrix whose 
eigenvalues $0$ (with multiplicity $m-1$) and $||\bx||^2$ have modulus less 
than $1$, so that $X$ is an $M$-matrix 
(see, e.g.,\cite{MO}). By the arguments in the last part of the 
proof of Theorem~\ref{t-3} we see that
$$\sum_{\bi\in Q_{m-1,m}}\det X(\bi,\bi)-m\det X=\text{tr} (I_m-X),$$
both sides being equal to $||\bx||^2$.
\end{remark}

Recall that if $C$ is an $n\times n$ matrix and $1\le k\le n$ one defines 
the $k$-th additive compound matrix $\De_{k}(C)$ as the 
$\binom{n}{k}\times\binom{n}{k}$ given by
$$\De_k(C)=\frac{d}{dt}(I_n+tC)^{(k)}\bigg|_{t=0}.$$
Equivalently, $\De_{k}(C)$ is the coefficient matrix of $t$ in the expansion
$(I_n+tC)^{(k)}=I_n+t\De_k(C)+\calO(t^2)$. As an operator, this is the same 
as the restriction to $\wedge^{k}\calH$ of the operator 
$$C\otimes I_n\otimes\cdots\otimes I_n+
I_n\otimes C\otimes I_n\otimes\cdots\otimes I_n
+\ldots +I_n\otimes I_n\otimes\cdots\otimes C$$
defined on $\otimes^{k}\calH$. 
It is well known that if $C$ has eigenvalues $\la_1,\ldots,\la_n$ then 
the eigenvalues of $\De_k(C)$ are $\la_{i_1}+\ldots+\la_{i_k}$, 
$1\le i_1<\ldots<i_k\le n$ (cf., e.g., \cite{MO}). 

\begin{proof}[Proof of Corollary~\ref{t-2}]
Using~\eqref{dec} and arguing as in the first part of the proof of 
Theorem~\ref{t-3} we arrive at
$$\frac{d}{dt}(I_n+tA)^{(k)}\bigg|_{t=0}
=U^{(k)}\frac{d}{dt}
\left(\text{diag}(1+tz_1,\ldots,1+tz_n)^{(k)}\right)\bigg|_{t=0}U^{*(k)}.$$
Since $A_n$ is upper triangular the diagonal elements of $(I_n+tA)^{(k)}$
corresponding to the $k\times k$ principal minors indexed by multi-indices
$\bi=(i_1,\ldots,i_k)\in Q_{k,n-1}$ have the form 
$\prod_{j=1}^{k}(1+tw_{i_j})$. The conclusion now follows from 
the obvious identities
$$\frac{d}{dt}\prod_{j=1}^{k}(1+tw_{i_j})\bigg|_{t=0}=\sum_{j=1}^{k}w_{i_j},
\quad \frac{d}{dt}\prod_{l=1}^{k}(1+tz_{s_l})=\sum_{l=1}^{k}z_{s_l}$$
combined with the arguments in the second 
part of the proof of Theorem~\ref{t-3}.
\end{proof}

Corollary~\ref{t-2} may alternatively be proved by using~\eqref{inv-m} with 
$\la=t$ and $\mu=1$ and differentiating with respect to $t$. Using this
alternative method we can actually prove Corollary~\ref{t-4}.

\begin{proof}[Proof of Corollary~\ref{t-4}]
Let $\la\in\bC\setminus\{0\}$ and $\mu\in\bC$. 
By~\eqref{inv-m} and Remark~\ref{la0} there exists an 
$\binom{n-1}{k}\times\binom{n}{k}$ row stochastic matrix $R_k$ such that
the relations
$$W_{k}^{[k]}(t\la,t\mu+1)^T=R_{k}Z_{k}^{[k]}(t\la,t\mu+1)^T\text{ and }
\bb^{[k]}=\ba^{[k]}R_{k}$$
hold for all $t\in\bC$. Thus
$$\frac{d^m}{dt^m}W_{k}^{[k]}(t\la,t\mu+1)^T\bigg|_{t=0}
=R_{k}\frac{d^m}{dt^m}Z_{k}^{[k]}(t\la,t\mu+1)^T\bigg|_{t=0}$$
for $m\in\{1,\ldots,k\}$, which is the same as~\eqref{gen-m}.
\end{proof}

\begin{remark}\label{rem-expl-distct}
The assumption that the points $z_i$, $1\le i\le n$, are {\em distinct} that 
was made at the beginning of \S \ref{s21} was not actually used in the 
above proofs. However, this condition is natural when viewing -- as we did -- 
these points as charged particles generating a resultant electrostatic force 
given by~\eqref{f}. Clearly, the latter can be appropriately rewritten even 
in the case of multiple points.
\end{remark}

\subsection{Proof of Theorems~\ref{t-5}--\ref{t-6}}\label{s32}

The inequality $h(\calW,\calZ)\leq h(\calW_e,\calZ)$ is trivially true since
$\calW\subseteq\calW_e$. To prove the interesting part of Theorem~\ref{t-5} 
we make use of a result on spectral variations of normal matrices 
which is usually referred to 
as the Bauer-Fike theorem (cf., e.g., \cite[Theorem VI.3.3]{Bh}).

\begin{lemma}\label{b-f}
Let $A$ and $B$ be $n\times n$ complex matrices. If $A$ is normal then
$$h(\Si(B),\Si(A))\le ||A-B||,$$
where $h$ denotes the relative Hausdorff distance from $\Si(B)$ to $\Si(A)$. 
\end{lemma}

Recall the orthonormal basis $(\bv_1,\ldots,\bv_n)$ of $\calH$ where the matrix
representation of the normal operator $A$ with $\Si(A)=\calZ$ is given 
by~\eqref{form}, that is, 
\begin{equation}\label{tri}
A=\begin{pmatrix}
A_{n} & \bx\\
\by^{*} & \ze
\end{pmatrix}.
\end{equation}
In the above decomposition $A_n$ is an $(n-1)\times(n-1)$ upper triangular 
matrix with $\Si(A_n)=\calW$, 
$\ze=\langle A\bv_n,\bv_n\rangle=\sum_{i=1}^{n}a_{i}z_{i}$, and 
$\bx,\by\in\bC^{n-1}$. Set $B=A_n\oplus \ze$. An easy computation now 
shows that
\begin{equation*}
\begin{split}
||A-B||^{2}&=||\bx||^2=||\by||^2
=\sum_{i=1}^{n-1}|\langle A\bv_{n},\bv_i\rangle|^2
=||A\bv_n||^2-|\langle A\bv_n,\bv_n\rangle|^2\\
&=\sum_{i=1}^{n}a_{i}|z_i|^2-|\ze|^2=\si_2(\bz;\ba)^2,
\end{split}
\end{equation*}
which combined with Lemma~\ref{b-f} yields the inequality stated in 
Theorem~\ref{t-5}.

To prove Theorem~\ref{t-6} note first that by Theorem~\ref{t-5} it 
remains to check that $h(\calZ,\calW_e)\le \si_2(\bz;\ba)$. The latter 
inequality is invariant under non-singular affine transformations of 
the complex plane and it is therefore enough to show that 
$h(\calZ,\calW_e)\le \si_2(\bz;\ba)$ whenever $\bz=(z_1,\ldots,z_n)\in\bR^n$. 
Since in this case $A$ is a Hermitian operator 
we get from~\eqref{tri} that $\bx:=(x_1,\ldots,x_{n-1})^t=\by$ and 
$A_n=\text{diag}(w_1,\ldots,w_{n-1})$. Suppose $z_i$ is 
such that 
\begin{equation}\label{ineq1}
\min_{1\le j\le n-1}|z_i-w_j|>\si_2(\bz;\ba)
\end{equation}
and let $\bu=(u_1,\ldots,u_n)^t\in\bC^n$ be a unit eigenvector of $A$ with 
eigenvalue $z_i$. Then by~\eqref{tri} we may rewrite the identity 
$A\bu=z_{i}\bu$ as
\begin{equation}\label{ineq2}
(\ze-z_i)u_{n}=-\sum_{j=1}^{n-1}\bar{x}_ju_j\text{ and }
(w_j-z_i)u_j=x_ju_n,\quad 1\le j\le n-1.
\end{equation}
Note that $\si_2(\bz;\ba)>0$ since $z_1,\ldots,z_n$ are assumed to be 
pairwise distinct.
Moreover, the fact that $||\bu||=1$ combined with~\eqref{ineq2} implies that
$u_n\neq 0$. From~\eqref{ineq1} and the last $n-1$ identities in~\eqref{ineq2}
we get
$$|u_j|\le \si_2(\bz;\ba)^{-1}|x_j||u_n|,\quad 1\le j\le n-1,$$
and then using the first relation in~\eqref{ineq2} we deduce that 
$$|\ze-z_i|\le \si_2(\bz;\ba)^{-1}||\bx||^2=\si_2(\bz;\ba).$$
This shows that $h(\{z_i\},\calW_e)\le \si_2(\bz;\ba)$, which completes 
the proof of Theorem~\ref{t-6}.

\subsection{Proof of Theorems~\ref{t-7}--\ref{t-8}}\label{s33}

It is clear that up to a rescaling we may assume that $\rho=1$ in both 
Theorem~\ref{t-7} 
and Theorem~\ref{t-8}. Below we shall use the Hilbert space of square 
summable complex sequences $\calH=l^2(\bN)$ 
with standard scalar product $\langle \cdot,\cdot \rangle$ and standard 
orthonormal basis $\{\bbe_i\}_{i\in\bN}$.

\begin{proof}[Proof of Theorem~\ref{t-7}]
We may assume wlog that $\{z_i\}_{i\in\bN}$ is an increasing sequence. 
The fact that $f$ has infinitely many zeros $\{w_j\}_{j\in\bN}$ that interlace
the poles $\{z_i\}_{i\in\bN}$ on the real axis follows by noticing that all 
zeros of $f$ are necessarily real and that $f$ is strictly decreasing on each
interval $(z_i,z_{i+1})$, $i\in\bN$. Alternatively, this may also be viewed as
a consequence of the operator theoretic interpretation given below combined 
with Hochstadt's interlacing theorem for the eigenvalues of 
rank one perturbations of compact selfadjoint Hilbert space operators 
\cite{H}. 

Let $A$ denote the unique bounded operator on $\calH$ satisfying 
$A\bbe_i=(1-z_i)\bbe_i$, $i\in\bN$. Note that $A$ is injective, 
selfadjoint and compact since the $z_i$ are real and $z_i\to 1$ as 
$i\to\infty$. By assumption one has $\sum_{i=1}^{\infty}a_i=1$ and so we may
define a unit vector $\bv\in\calH$ in the following way
\begin{equation}\label{i-v}
\bv=\sum_{i=1}^{\infty}\sqrt{a_i}\bbe_i.
\end{equation}
Set $\calK=\bv^{\perp}$ and let $P$ be the orthoprojection on $\calK$. Recall
from \cite[\S IV.6]{Ka2} that the Weinstein-Aronszajn determinant of the 
second kind associated with $A$ and $P$ is given by 
$\langle (A-zI)^{-1}\bv,\bv\rangle$. A simple computation shows that
\begin{equation}\label{wa1}
\langle (A-zI)^{-1}\bv,\bv\rangle=f(1-z),\quad z\in\bC\setminus 
\{1-z_i\}_{i\in\bN}.
\end{equation}
Let us now consider the compression $A'$ of $A$ to $\calK$, that is, the 
operator on $\calK$ given by $A'=PAP|_{\calK}$. Note that $A'$ is itself a 
compact selfadjoint operator. We claim that the point spectrum of $A'$ consists
precisely of the zeros of $f(1-z)$.

\begin{lemma}\label{wa2}
The set of eigenvalues of $A'$ is $\{1-w_j\}_{j\in\bN}$.
\end{lemma}

\begin{proof}
If $\bu\in\calK$ is an eigenvector of $A'$ with eigenvalue $\ze\in\bC$ 
(actually $\bR$) then
\begin{equation}\label{wa3}
(A-\ze I)\bu=(I-P)A\bu +(A'-\ze I)\bu=\langle A\bu,\bv\rangle\bv.
\end{equation}
From~\eqref{wa3} we see that $\langle A\bu,\bv\rangle\neq 0$ since otherwise 
$\bu$ would be an eigenvector of $A$, say $\bu=\bbe_i$ for some $i\in\bN$, 
which would then lead to the contradiction 
$\langle A\bu,\bv\rangle=\langle A\bbe_i,\bv\rangle=(1-z_i)\sqrt{a_i}\neq 0$. 
By~\eqref{wa1} and~\eqref{wa3} we get
$$f(1-\ze)=\langle (A-\ze I)^{-1}\bv,\bv\rangle
=\langle A\bu,\bv\rangle^{-1}\langle \bu,\bv\rangle=0$$
so that $1-\ze=w_j$ for some $j\in\bN$. Conversely, let $w$ be such that 
$f(1-w)=0$, that is, $\langle (A-wI)^{-1}\bv,\bv\rangle=0$, and set 
$\bu=(A-wI)^{-1}\bv$. Then clearly $\bu\neq 0$ and one has
$$\bv=(A-wI)\bu=\langle A\bu,\bv\rangle\bv+(A'-wI)\bu.$$
Since $(A'-wI)\bu\in\calK$ and $\bv\in\calK^{\perp}$ it follows that
$$(A'-wI)\bu=(1-\langle A\bu,\bv\rangle)\bv=0.$$
Thus $\bu$ is an eigenvector of $A'$ with eigenvalue $w$.
\end{proof}

It is not difficult to see that $A'$ is injective. Indeed, if 
$0\neq\bu\in\calK$ is such that
$A'\bu=0$ then from $PA\bu=0$ and the fact that $A$ is injective we get
$A\bu=\al\bv$ for some $\al\in\bC\setminus\{0\}$. Hence
$$0=\langle \bu,\bv\rangle=\al\langle A^{-1}\bv,\bv\rangle=\al f(1),$$
which is a contradiction since the assumptions of the theorem (more 
specifically, condition~\eqref{i-hyp} with $\rho=1$) imply that $f(1)>0$. 

By the spectral theorem for compact selfadjoint operators 
\cite[Corollary X.3.3.5]{DS} there exists a complete orthonormal 
system $\{\bbf_j\}_{j\in\bN}$ in $\calK$ consisting of eigenvectors for $A'$.
Set $s_{jk}=|\langle \bbf_j,\bbe_k\rangle|^2$ for $j,k\in\bN$. Then
$$1-w_j=\langle A'\bbf_j,\bbf_j\rangle=\langle A\bbf_j,\bbf_j\rangle
=\sum_{k=1}^{\infty}s_{jk}(1-z_k).$$
Both $\{\bv\}\cup\{\bbf_j\}_{j\in\bN}$ and $\{\bbe_k\}_{j\in\bN}$ are 
complete orthonormal systems in $\calH$, so that 
\begin{equation*}
\begin{split}
&\sum_{k=1}^{\infty}s_{jk}=\sum_{k=1}^{\infty}|\langle \bbf_j,\bbe_k\rangle|^2
=||\bbf_j||^2=1,\quad j\in\bN,\\
&\sum_{j=1}^{\infty}s_{jk}=\sum_{j=1}^{\infty}|\langle \bbf_j,\bbe_k\rangle|^2
=||\bbe_k||^2-|\langle \bv,\bbe_k\rangle|^2=1-a_k,\quad k\in\bN,
\end{split}
\end{equation*}
which completes the proof of the theorem.
\end{proof}

\begin{proof}[Proof of Theorem~\ref{t-8}]
Let $\{z_i\}_{i\in\bN}$ be a sequence of distinct complex numbers 
satisfying~\eqref{i-hyp} with $\rho=1$ and let $A$ be the injective compact 
normal operator on $\calH$ defined by $A\bbe_i=(1-z_i)\bbe_i$, $i\in\bN$. As
before we denote by $A'$ the compression $PAP|_{\calK}$ of $A$, where 
$\calK=\bv^{\perp}$ and $\bv$ is given by~\eqref{i-v}. Clearly, $A'$ is a 
compact operator on $\calK$. Moreover, since $\Re f(1)>0$ the arguments in 
the proof of Theorem~\ref{t-7} show that $A'$ is injective. Note that 
Lemma~\ref{wa2} remains valid in the present setting and that by assumption 
$f$ has a discrete set of zeros $\{w_j\}_{j\in\bN}$. It follows that the 
point spectrum of $A'$ is $\{1-w_j\}_{j\in\bN}$. 
Denote by $\calK_{A'}$ the closed linear hull of the root subspaces of $A'$ 
and let $\widehat{A'}$ be the operator on $\calK_{A'}$ induced by $A'$. Recall 
the operator version of Schur's lemma given in \cite[Lemma 4.1]{GK}.

\begin{lemma}\label{sch}
With the above notations there exists an orthonormal basis 
$\{\bbf_j\}_{j\in\bN}$ of $\calK_{A'}$ in which the matrix of the operator 
$\widehat{A'}$ has triangular form and 
$\langle \widehat{A'}\bbf_j,\bbf_j\rangle=1-w_j$, $j\in\bN$. 
\end{lemma}

Set $s_{jk}=|\langle \bbf_j,\bbe_k\rangle|^2$ for $j,k\in\bN$. By 
Lemma~\ref{sch} we have
$$1-w_j=\langle \widehat{A'}\bbf_j,\bbf_j\rangle=\langle A\bbf_j,\bbf_j\rangle
=\sum_{k=1}^{\infty}s_{jk}(1-z_k).$$
It is clear that
\begin{equation}\label{bes-1}
\sum_{k=1}^{\infty}s_{jk}=\sum_{k=1}^{\infty}|\langle \bbf_j,\bbe_k\rangle|^2
=||\bbf_j||^2=1,\quad j\in\bN.
\end{equation}
On the other hand $\{\bv\}\cup\{\bbf_j\}_{j\in\bN}$ is an orthonormal system 
in $\calH$ so that by Bessel's inequality we get
\begin{equation}\label{bes-2}
\sum_{j=1}^{\infty}s_{jk}=\sum_{j=1}^{\infty}|\langle \bbf_j,\bbe_k\rangle|^2
\le ||\bbe_k||^2-|\langle \bv,\bbe_k\rangle|^2=1-a_k,\quad k\in\bN.
\end{equation}
The inequality stated in the theorem by means of nonnegative convex 
functions is an immediate consequence of~\eqref{bes-1}--\eqref{bes-2}. 
\end{proof}

\section{A hierarchy of de Bruijn-Springer relations and inertia 
laws}\label{s-mom} 

Theorems~\ref{t-1}--\ref{t-3} suggest even deeper physical connections 
between the set of equilibrium points and the given set of charged particles 
in the case of finite planar systems, such as a whole new hierarchy of 
inertia laws for planar solid bodies associated with these sets. In order to
make a precise statement 
let us identify the complex plane $\bC$ with 
the $x_1x_2$-plane in $\bR^3$ equipped with the standard coordinate system 
$(O;x_1,x_2,x_3)$. For an integer $k\ge 2$ we denote by $\Xi_k$ the standard 
simplex
$$\Xi_k=\left\{(t_1,\ldots,t_k)\in [0,1]^k\mid t_1+\ldots +t_k
=1\right\}.$$

\begin{definition}
Let $L$ be a line in $\bR^3$. Given a convex $k$-gon $\calK$ in $\bC$ with 
vertices $v_1,\ldots,v_k$ and $\al\ge 1$ we define the $\al$-{\em moment 
of} $\calK$ {\em with respect to} $L$ to be 
$$I_{\al,L}(\calK)=\int\cdots\int_{\Xi_k}
d\!\left(t_{1}v_1+\ldots +t_{k}v_k,L\right)^{\al}dt_1\cdots dt_k,$$
where $d(\cdot,L)$ denotes the distance function associated with $L$.
\end{definition}

\begin{remark}
If $\calK$ is viewed as a planar solid body with uniform mass distribution then
$I_{\al,L}(\calK)$ may be interpreted as a higher moment of inertia (or 
generalized angular momentum) of $\calK$ with respect to the rotation axis $L$.
\end{remark}

Let $\Si_m$ denote the symmetric group on $m$ elements. It is not difficult to 
see that if $\calK$ is as above and $u_1,\ldots,u_m\in\bC$ are such that 
$\calK=\co(u_1,\ldots,u_m)$ then 
\begin{equation}\label{invar-m}
\int\cdots\int_{\Xi_m}d(t_{\pi(1)}u_1+\ldots +t_{\pi(m)}u_m,L)^{\al}
dt_1\cdots dt_m=I_{\al,L}(\calK)
\end{equation}
for any $\pi\in \Si_m$, so that $I_{\al,L}(\calK)$ depends only on the extreme 
points of $\calK$. 

The following conjecture is a powerful generalization of 
Theorems~\ref{t-1}--\ref{t-3} and suggests a ``hierarchy'' of weighted 
de Bruijn-Springer relations and inertia laws for planar polygons constructed 
from the equilibrium points and the charged particles. 

\begin{conjecture}\label{gen-mom}
Let $z_i$, $1\le i\le n$, be distinct points in $\bC$ and $a_i>0$, 
$1\le i\le n$, be such that $\sum_{i=1}^{n}a_i=1$, where $n\ge 2$. Denote by 
$w_j$, $1\le j\le n-1$, the zeros of 
$$f(z)=\sum_{i=1}^{n}\frac{a_i}{z-z_i}$$
and let $\Phi:\bC\to \bR$ be an arbitrary convex function. Then for any 
$k\in\{1,\ldots,n-1\}$, $1\le m\le k$, and $(t_1,\ldots,t_k)\in\bC^k$ one has
\begin{equation}\label{dbs-h}
\begin{split}
\sum_{1\le r_1<\ldots<r_k\le n-1}&\,\,\sum_{\pi\in \Si_k}
\Phi\left(\Pi_{k,m}(t_{\pi(1)}w_{r_1},\ldots,t_{\pi(k)}w_{r_k})\right)\\
\le \sum_{1\le s_1<\ldots <s_k\le n}&\left(1-\sum_{i=1}^{k}a_{s_i}\right)
\sum_{\pi\in \Si_k}
\Phi\left(\Pi_{k,m}(t_{\pi(1)}z_{s_1},\ldots,t_{\pi(k)}z_{s_k})\right),
\end{split}
\end{equation}
where $\Pi_{k,m}$ is the $m$th elementary symmetric function on $k$ symbols.
In particular, for any line $L$ in $\bR^3$ and $\al\ge 1$ one has
\begin{equation}\label{ang-m}
\begin{split}
\sum_{1\le r_1<\ldots<r_k\le n-1}
& I_{\al,L}\left(\tco(w_{r_1},\ldots,w_{r_k})\right)\\
& \le \sum_{1\le s_1<\ldots <s_k\le n}\left(1-\sum_{i=1}^{k}a_{s_i}\right)
I_{\al,L}\left(\tco(z_{s_1},\ldots,z_{s_k})\right).
\end{split}
\end{equation}
\end{conjecture}

The fact that~\eqref{dbs-h} implies~\eqref{ang-m} follows from~\eqref{invar-m}
by taking $m=1$ and integra\-ting over the simplex $\Xi_k$. Arguing as in the 
proof of Corollary~\ref{t-4} one can see that it is actually enough to 
prove~\eqref{dbs-h} for $m=k$ and $k\in\{1,\ldots,n-1\}$. Note also that in 
the case of equal charges Conjecture~\ref{gen-mom} would provide a whole new 
series of inequalities for averages of generalized inertia moments 
involving 
the critical points and the zeros of an arbitrary complex polynomial. These 
would substantially improve both Corollary~\ref{cor-mom} and the 
inequalities conjectured in \cite[Remark 2.3.9]{RS}.  

\begin{remark}
In view of the proofs of \cite[Theorem 4.7]{Ma2} and Theorem~\ref{t-3} above,
a natural way of attacking Conjecture~\ref{gen-mom} would be to interpret
the arguments of $\Phi$ in~\eqref{dbs-h} as eigenvalues of some appropriately
defined multilinear power ``$\wedge^k A$'' of a normal matrix $A$ and one of 
its principal submatrices, where ``$\wedge$'' is a suitable multilinear 
functional on $\otimes^k \calH$. However, defining such a functional appears 
to be a non-trivial question. Note for instance that an obvious polarization
of the usual antisymmetric tensor power, i.e., setting 
$$A\wedge B=\frac{1}{2}\left[(A+B)^{(2)}-A^{(2)}-B^{(2)}\right]$$ 
for square matrices $A,B$ of the same order, cannot lead to the desired
conclusion since to begin with one does not get the correct dimensions (these
should agree with the number of $k$-tuples in the two sides of~\eqref{dbs-h}).
\end{remark}

\section{Zeros of Borel series with positive $l^1$-coefficients}\label{s42}

Theorem~\ref{t-8} applies essentially to any bounded discrete set of 
positively charged particles with finite total charge for which equilibrium 
points are known to exist. We conjecture that if the particles accumulate at 
a point on the boundary of a disk containing all the particles then there 
exist infinitely many equilibrium points.

\begin{conjecture}\label{con-2}
Let $\{a_i\}_{i\in\bN}$ be a sequence of positive numbers whose sum is finite
and $\{z_i\}_{i\in\bN}$ a sequence of distinct complex numbers 
satisfying
$$|z_i|<\rho\text{ for }i\in\bN\text{ and }\lim_{i\to\infty}z_i
=\xi\text{ with }|\xi|=\rho$$
for some fixed positive number $\rho$. Then the function
\begin{equation}\label{i-ff}
f(z)=\sum_{i=1}^{\infty}\frac{a_i}{z-z_i}
\end{equation}
has an infinite set of zeros in the disk $\bD(\rho)=\{z\in\bC: |z|<\rho\}$.
\end{conjecture} 

Note that the assumptions of Conjecture~\ref{con-2} imply 
that condition~\eqref{i-hyp} holds and thus $f$ is meromorphic in 
$\bD(\rho)$. Functions representable as a series of the form~\eqref{i-ff} are
known as Borel series \cite{W} and so Conjecture~\ref{con-2} may be rephrased
as follows: {\em every Borel series with positive $l^1$-coefficients whose 
poles
lie in an open disk and accumulate at a point on the boundary of this disk 
has infinitely many zeros}. 

Conjecture~\ref{con-2} may also be translated in operator theoretic terms by 
using the methods employed in the proofs of Theorems~\ref{t-7}--\ref{t-8}. 
Note first that up to a similarity transformation of the complex plane we
may assume that $\xi=\rho=1$. Then one can show that 
the above conjecture is equivalent to the following statement.

\begin{conjecture}\label{new-c2}
Let $\{z_n\}_{n\in\bN}$ be a sequence of distinct complex numbers in the open 
unit disk converging to $1$ and denote by $A$ the (bounded) normal operator on 
a separable infinite-dimensional complex Hilbert space $\calH$ satisfying 
$A\bbe_n=z_{n}\bbe_n$, $n\in \bN$, where $\{\bbe_n\}_{n\in\bN}$ is an 
orthonormal basis of $\calH$. If $P$ is a bounded operator on $\calH$ such 
that $I-P$ is a rank one orthoprojection then the compression 
$PAP|_{P\calH}$ has infinite point spectrum.
\end{conjecture}

\begin{remark}\label{more-gen}
The hypotheses in Conjecture~\ref{new-c2} imply that the normal operator 
$I-A$ is compact, which may facilitate its study. However, it seems likely 
that Conjecture~\ref{con-2} should actually 
hold even under the less restrictive conditions $|z_i|<\rho$ for all 
$i\in\bN$ and $\lim_{i\to\infty}|z_i|=\rho$. In other words, it is reasonable 
to believe that the convergence condition in Conjecture~\ref{new-c2} may be 
relaxed 
as follows: $|z_n|<1$ for all $n\in\bN$ and $\lim_{n\to\infty}|z_n|=1$.
\end{remark}

\section{Operator versions of the Clunie-Eremenko-Rossi 
conjecture}\label{s43}

Various results concerning the existence of equilibrium 
points for Newtonian and logarithmic potentials have been obtained in
e.g.~\cite{CER}, \cite{ELR} and \cite{LR}. It should be emphasized though 
that these deal almost exclusively with {\em unbounded} discrete charge 
configurations. Note that Conjecture~\ref{con-2} 
may actually be viewed as a natural analog for 
{\em bounded} discrete charge configurations of the following well-known 
conjecture
of Clunie-Eremenko-Rossi, see, e.g., \cite[Conjecture 2.7]{CER} and 
\cite[Conjecture 1.1]{LR}.

\begin{conjecture}[Clunie-Eremenko-Rossi conjecture]\label{con-CER}
Let $\{a_i\}_{i\in\bN}$ be a sequence of positive numbers and 
$\{z_i\}_{i\in\bN}$ a sequence of distinct complex numbers such that
$$z_i\to\infty\text{ as }i\to\infty\text{ and }
\sum_{z_i\neq 0}\frac{a_i}{|z_i|}<\infty.$$
Then the meromorphic function 
$$f(z)=\sum_{i=1}^{\infty}\frac{a_i}{z-z_i}$$
has infinitely many zeros.
\end{conjecture}

A catchy albeit somewhat loose rephrasing of this conjecture is as follows:
{\em every flat universe has infinitely many resting points}. 

Let us use the framework of \S \ref{s33} to give an operator theoretic
interpretation and generalization of Conjecture~\ref{con-CER}. Note first that 
up to a translation of the variable $z$ and a rescaling of the coefficients
$a_n$ we may assume wlog that all $z_n$ are non-zero and satisfy the 
normalization condition
\begin{equation}\label{i-norm}
\sum_{i=1}^{\infty}\frac{a_i}{|z_i|}=1.
\end{equation}
As in \S \ref{s33} we use the Hilbert space $\calH=l^2(\bN)$ with standard 
inner 
product $\langle\cdot,\cdot\rangle$ and standard complete orthonormal system
$\{\bbe_n\}_{n\in\bN}$. Set $z_n=r_n e^{i\te_n}$, $n\in\bN$, and let $A$ and
$B$ be the bounded normal operators on $\calH$ satisfying
\begin{equation}\label{ab1}
A\bbe_n=e^{i\te_n}\bbe_n,\quad B_n=r_{n}^{-1}\bbe_n,\quad n\in\bN.
\end{equation}
Clearly, $A$ is unitary and $B$ is compact. Condition~\eqref{i-norm} allows us 
to define a unit vector $\bv\in\calH$ by setting
$$\bv=\sum_{n=1}^{\infty}\sqrt{\frac{a_n}{|z_n|}}\,\bbe_n$$
and one can easily verify that
\begin{equation}\label{ab2}
\langle (A-zB)^{-1}\bv,\bv\rangle=-f(z),\quad 
z\in\bC\setminus \{z_n\}_{n\in\bN}.
\end{equation}
Set $\calK=\bv^{\perp}$ and let $P$ be the orthoprojection on $\calK$. We use
the notation $A'$ and $B'$ for the compressions to $\calK$ of $A$ and $B$, 
respectively, that is,
$$A'=PAP|_{\calK},\quad B'=P^BP|_{\calK}.$$
By definition, the set of generalized eigenvalues of $A'$ with respect to
$B'$ is
$$\Si_{_{B'}}(A')=\{z\in\bC\mid N(A'-zB')\neq \{0\}\},$$
where $N(T)$ denotes as usual the null space of an operator $T$. Let $\calZ(f)$
be the zero set of $f$. 

\begin{lemma}\label{i-eig}
In the above notation one has $\calZ(f)=\Si_{_{B'}}(A')$.
\end{lemma}

\begin{proof}
We make an appropriate adaption of the argument used in the proof of 
Lemma~\ref{wa2}. Using the fact that the $z_n$ are distinct 
one can check that if $0\neq \bu\in\calH$ then $\bu\in N(A-\ze B)$ if and 
only if $\ze=z_n$ and $\bu=\bbe_n$ for some $n\in\bN$. However
$$\langle\bbe_n,\bv\rangle=\sqrt{\frac{a_n}{|z_n|}}\neq 0,$$
which shows that $N(A-\ze B)\cap \calK=\{0\}$ for $\ze\in\bC$. Let now 
$\bu\in\calK$ be such that $(A'-\ze B')\bu=0$ for some $\ze\in\bC$. From the
identity
$$(A-\ze B)\bu=(I-P)(A-\ze B)\bu+(A'-\ze B')\bu
=\langle (A-\ze B)\bu,\bv\rangle\bv$$
and $\langle (A-\ze B)\bu,\bv\rangle\neq 0$ it follows that 
$(A-\ze B)^{-1}\bv=\langle (A-\ze B)\bu,\bv\rangle^{-1}\bu$. 
By~\eqref{ab2} one gets
$$f(\ze)=-\langle (A-\ze B)\bu,\bv\rangle^{-1}\langle \bu,\bv\rangle=0,$$
so that $\ze\in\calZ(f)$. Conversely, let $w\in\calZ(f)$. Then 
$w\notin\{z_n\}_{n\in\bN}$ and so $A-wB$ is invertible by~\eqref{ab1}. The 
vector $\bu:=(A-wB)^{-1}\bv$ is clearly non-zero and one has 
$\langle \bu,\bv\rangle=f(w)=0$ hence $\bu\in\calK$. Thus
$$\bv=(A-wB)\bu=\langle (A-wB)\bu,\bv\rangle\bv+(A'-wB')\bu,$$
which shows that $(A'-wB')\bu\in\bv^{\perp}\cap\{\bv\}=\{0\}$ and therefore
$w\in \Si_{_{B'}}(A')$.
\end{proof}

Recall that a vector with all non-zero coordinates in the standard 
basis $\{\bbe_n\}_{n\in\bN}$ of $\calH$ is called totally non-zero. The 
above discussion shows that Clunie-Eremenko-Rossi conjecture 
(Conjecture~\ref{con-CER}) is in fact 
equivalent to the following statement.

\begin{conjecture}\label{con-5}
Let $\{\al_n\}_{n\in\bN}$ and $\{\be_n\}_{n\in\bN}$ be a sequence of 
unimodular complex numbers and a sequence of positive numbers converging to 
$0$, respectively. Assume that $\al_{n}\be_{n}^{-1}$ are distinct and 
define a unitary operator $A$ and a compact selfadjoint 
operator $B$ on $\calH$ by setting
$$A\bbe_n=\al_n \bbe_n,\quad B\bbe_n=\be_n \bbe_n,\quad n\in\bN.$$
Let $\bv$ be a totally non-zero unit vector in $\calH$, denote by $P$ the 
orthoprojection on $\calK=\bv^{\perp}$ and let $A'=PAP|_{\calK}$ and 
$B'=PBP|_{\calK}$. Then $\left|\Si_{_{B'}}(A')\right|=\infty$, i.e., 
there exist infinitely many $z\in\bC$ such that $N(A'-zB')\neq\{0\}$.
\end{conjecture}

If the vector $\bv$ in Conjecture~\ref{con-5} is represented in the standard 
orthonormal basis of $\calH$ as 
$\bv=\sum_{n=1}^{\infty}v_n\bbe_n$ then Keldysh's value distribution theorem 
for meromorphic functions implies that 
Conjecture~\ref{con-5} is true provided that 
$\sum_{n=1}^{\infty}\be_{n}^{-1}|v_n|^2<\infty$ (see, e.g., \cite{CER} and 
\cite[Ch.~5, Theorem 6.2]{GO}. 

Conjecture~\ref{con-5} suggests the following more general question: let $A$ 
be an injective compact normal Hilbert space operator 
and let $B$ and $C$ be finite-rank operators. Is it true that 
$\left|\Si_{_{I+C}}(A+B)\right|=\infty$, that is, are there 
infinitely many $z\in\bC$ such that $N(A+B-z(I+C))\neq\{0\}$? It is not 
difficult to see that an affirmative answer to this question would actually 
imply the validity of the operator theoretic version of the 
Clunie-Eremenko-Rossi 
conjecture (Conjecture~\ref{con-5}). Indeed, in the notation of 
Conjecture~\ref{con-5} we may write
\begin{equation*}
\begin{split}
A'-zB'=-zA\big\{A^{-1}B-A^{-1}[&PB(I-P)-(I-P)B]\\
&-z^{-1}[I-A^{-1}(PA(I-P)+(I-P)A)]\big\}.
\end{split}
\end{equation*}
It remains to note that under the assumptions of 
Conjecture~\ref{con-5} the operator 
$A^{-1}B$ is injective, compact and normal while both 
$A^{-1}[PB(I-P)-(I-P)B]$
and $A^{-1}[PA(I-P)+(I-P)A]$ have rank at most two. However, the answer to 
the question raised above is negative, as 
one can see from the following example.

\begin{ex} 
Let $J$ be the simple dissipative 
Volterra operator on $L^{2}(0,1)$ defined by
$$(Jf)(x)=2i\int_{0}^{x}f(t)dt$$
and set $A=J+J^{*}$, $B=J-J^{*}$ and $C=0$. Then $A$ is 
injective compact selfadjoint, $B$ is a one-dimensional operator while 
the point spectrum of $A+B=2J$ is empty (cf.~\cite[p.~187]{GK}), so that
in this case $A+B-z(I+C)$ is injective for all $z\in\bC$.
\end{ex}

Let us finally mention that the existence of equilibrium points for some  
discrete potentials generated by complex charges 
was established in e.g.~\cite{LR}. We propose the following analog of 
Conjecture~\ref{con-CER} for certain potentials of this type.

\begin{conjecture}\label{con-compl}
Let $\{a_n\}_{n\in\bN}$ and $\{z_n\}_{n\in\bN}$ be sequences of complex numbers
such that the $z_n$ are distinct and satisfy
$$a_{n}\bar{z}_n>0\text{ for }n\in\bN,\quad z_n\to\infty\text{ as }n\to\infty,
\quad \sum_{n=1}^{\infty}\left|\frac{a_n}{z_n}\right|<\infty.$$
Then the meromorphic function
$$f(z)=\sum_{n=1}^{\infty}\frac{a_n}{z-z_n}$$
has infinitely many zeros.
\end{conjecture}

As in the case of Conjectures~\ref{con-2} and~\ref{con-CER}, one can 
reduce the conjecture stated above to a purely operator theoretic question. 
Indeed, in the notations of Conjecture~\ref{con-compl} let $A$ denote the 
injective compact normal operator on $\calH$ satisfying 
$A\bbe_n=z_{n}^{-1}\bbe_n$, $n\in\bN$. Define a vector
$$\bv=\sum_{n=1}^{\infty}\sqrt{\left|\frac{a_n}{z_n}\right|}\,\bbe_n\in\calH,$$
let $P$ the orthoprojection on $\bv^{\perp}$ and set 
$A'=PAP|_{\bv^{\perp}}$. Then
$$\langle (A-z^{-1}I)^{-1}\bv,\bv\rangle=zf(z),\quad 
z\in\bC^{*}\setminus \{z_n\}_{n\in\bN}.$$
Arguing as in the proof of Lemma~\ref{i-eig} one can show that 
$\ze^{-1}\in\calZ(f)$ for any $\ze\in\Si_{p}(A')\setminus \{0\}$, where 
$\calZ(f)$ is the zero set of $f$ and $\Si_{p}(A')$ denotes the point spectrum 
of $A'$. For Conjecture~\ref{con-compl} it would therefore be enough to show 
that $\left|\Si_{p}(A')\right|=\infty$.

\section*{Acknowledgements}

The author is grateful to Mihai Putinar, Nikolai Nikolski, Alexandre Eremenko, 
Donald Sarason and Boris Shapiro for their interest in this work and 
valuable remarks as well as to Semen Malamud for making available his 
inspiring paper~\cite{Ma2}. Thanks are also due to the anonymous referee 
for numerous constructive suggestions and for pointing out several useful 
references.

\end{document}